\documentclass[11pt]{article}
\usepackage{psfrag,graphicx}
\usepackage{amsfonts,amssymb}
\usepackage{epsf,latexsym}
\input epsf
\begin{document}
\newfont{\Bb}{msbm10}
\def\ni{\noindent}
\def\be{\begin{enumerate}}
\def\ee{\end{enumerate}}
\def\ec{\end{center}}
\def\bc{\begin{center}}
\def\-1{{^{-1}}}
\def\th{{\theta}}
\def\f{{\phi}}
\def\l{{\lambda}}
\def\w{{\wedge}}
\def\x{{\times}}
\def\ff{{\tilde{f}}}
\def\FF{{\tilde{F}}}
\def\Proof{{\bf Proof}\hspace{.3in}}
\def\a{{\alpha}}
\def\b{{\beta}}
\def\g{{\gamma}}
\def\G{{\Gamma}}
\def\o{{\omega}}
\def\i{{\iota}}
\def\t{{\theta}}
\def\bolk{{\O^k_V/h\O^k_V(\log D)}}
\def\olk{{\frac{\O^k_V}{h\O^k_V(\log D)}}}
\def\olo{{\frac{\O^0_V}{h\O^0_V(\log D)}}}
\def\olu{{\frac{\O^1_V}{h\O^1_V(\log D)}}}
\def\old{{\frac{\O^2_V}{h\O^2_V(\log D)}}}
\def\olm{{\frac{\O^{p-1}_V}{h\O^{p-1}_V(\log D)}}}
\def\X{{\cal X}}
\def\O{{\Omega}}
\def\ch{{\check{\Omega}}}
\def\chch{{\check{\cal H}}}
\def\chh{{\check{H}}}
\def\k{{\kappa}}
\def\p{{\partial}}
\def\CC{\mbox{\Bb C}}
\def\Ct{{\CC^p}}
\def\CT{{\CC^p\times\CC^d}}
\def\Cs{{\CC^n}}
\def\CS{{\CC^n\times \CC^d}}
\def\RR{\mbox{\Bb R}}
\def\HH{\mbox{\Bb H}}
\def\H{{\cal H}}
\def\PP{\mbox{\Bb P}}
\def\ZZ{\mbox{\Bb Z}}
\def\eop{\hspace*{\fill}$\Box$ \vskip \baselineskip}
\def\Der{\mbox{Der}}
\def\Ext{\mbox{Ext}}
\def\Ker{\mbox{Ker}}
\def\Coker{\mbox{Coker}}
\def\Im{\mbox{Im}}
\def\adj{\mbox{adj}}
\def\alg{\mbox{alg}}
\def\depth{\mbox{depth}}
\def\supp{\mbox{supp}}
\def\Tor{\mbox{Tor}}
\def\codim{\mbox{codim}}
\def\Hom{\mbox{Hom}}
\def\DD{{\Der(\log D)}}
\def\Dh{{\Der(\log h)}}
\def\pd{\mbox{pd}}
\def\cod{${\cal A}_e\codim$}
\def\D{{\cal D}}
\def\EE{{\cal E}}
\def\R{{\cal R}}
\def\F{{\cal F}}
\def\L{{\cal L}}
\def\B{{\cal B}}
\def\C{{\cal C}}
\def\A{{\cal A}}
\def\I{{\cal I}}
\def\T{{\cal T}}
\def\OO{{\cal O}}
\def\Ot{{{\cal O}_{\CC^p}}}
\def\OT{{{\cal O}_{\CC^p\times\CC^d}}}
\def\Os{{{\cal O}_{\CC^n}}}
\def\OS{{{\cal O}_{\CC^n\times\CC^d}}}
\def\iff{{\Leftrightarrow}}
\def\implies{{\Rightarrow}}
\def\hook{{\hookrightarrow}}
\def\to{{\ \rightarrow\ }}
\def\too{{\ \longrightarrow\ }}
\def\b{{\bullet}}
\def\u{\underline}
\def\dim{{\mbox{dim}}}
\newtheorem{theorem}{Theorem}[section]
\newtheorem{lemma}[theorem]{Lemma}
\newtheorem{theorem/definition}[theorem]{Theorem/Definition}
\newtheorem{proposition}[theorem]{Proposition}
\newtheorem{example}[theorem]{Example}
\newtheorem{remark}[theorem]{Remark}
\newtheorem{corollary}[theorem]{Corollary}
\newtheorem{definition}[theorem]{Definition}
\newtheorem{conjecture}[theorem]{Conjecture}
\title{Differential forms on free and almost free divisors}
\author{David Mond,\\
University of Warwick}
\maketitle
\hfill To Egbert Brieskorn\\

\begin{abstract}
We introduce a variant of the usual K\"{a}hler forms on singular free divisors,
and show that they enjoy the same depth properties as K\"{a}hler forms
on isolated hypersurface singularities. 
Using these forms it is possible to describe analytically the 
vanishing cohomology in families of free divisors, in precise analogy with
the classical description for the Milnor fibration of an ICIS, 
due to Brieskorn and Greuel. This applies in particular to the family
$\{D(f_\lambda)\}_{\lambda\in \Lambda}$ of discriminants of a versal deformation
$\{f_\lambda\}_{\lambda\in\Lambda}$ of a singularity of a mapping. 
\end{abstract}
\bc
\section{Introduction}
\label{introduction}
\ec
The hypersurface $D$ in the complex manifold $X$ is
a {\em free divisor} if the $\OO_X$-module
$\Der(\log D)$ of germs of vector fields on $X$ which are tangent to $D$,
is locally free( see \cite{saito}). 
Examples of free divisors include smooth hypersurfaces, 
normal crossing divisors, reflection arrangements (\cite{orlik and terao})
and discriminants of 
right-left 
stable maps with source dimension not less than
target dimension (\cite{looijenga}, 6.13). 

In many contexts, one is give a divisor $D_0$ which is free outside $0$,
and a family $D\to \CC$ with special fibre $D_0$, 
in which the general fibre is free. Typically, $D_t$ grows homology 
classes which vanish when $t$ returns to $0$.

A familiar example of this phenomenon is shown in the following pair
of diagrams. In the first, we see four planes passing through $0$ in 
3-space, of which each three are in general position.

\vskip 10pt
\epsfxsize 9cm
\centerline{\epsfbox{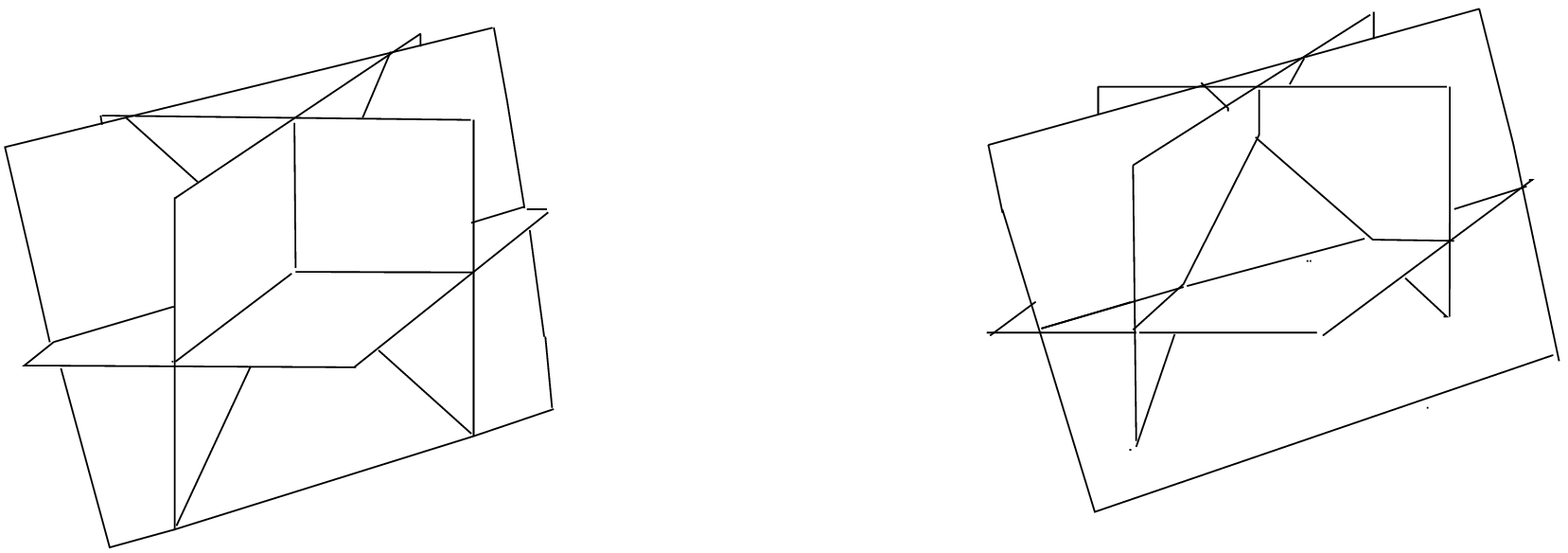}}
\vskip 10pt

Their union, $D_0$, is a normal crossing divisor, and
therefore free, outside $0$. It is not free at $0$; this can be
checked by a calculation, but it also follows from 2.2 below. In the
second picture, one of the planes has been shifted off $0$, and the union,
$D_t$, now has a non-trivial second homology class, 
carried by the tetrahedron one can see in the centre of the picture
(in fact the inclusion of the real $D_0$ and $D_t$ in their complexification,
induces a homotopy-equivalence, so these are ``good real pictures'').
Since $D_t$ is now a normal crossing divisor (everywhere), it is free.
One can think of $D_t$ as a kind of singular Milnor fibre of $D_0$ ---
we have not smoothed $D_0$ (as in ``classical'' singularity theory);
instead we have freed it.

Another large class of examples: if $f_0:(\CC^n,0)\to(\CC^p,0)$ (with
$n\geq p$)
is a germ of analytic map stable outside $0$ 
(equivalently, having ``finite
${\cal A}$-codimension'') then the discriminant $D(f_0)$ is free outside $0$;
if $f_t$ is a stable perturbation of $f_0$ then 
$D(f_t)$ is a free divisor, which carries homology classes which vanish when 
$t$ returns to $0$. 

In \cite{damonafd}
Damon introduced a definition of {\it almost free divisor}, and a suitable
category of deformations, which encompass
these examples and many others. Using Morse-theoretic techniques introduced in 
\cite{dm}, he
is able to calculate the rank of the vanishing homology. In this paper we 
describe the vanishing cohomology in the singular Milnor fibres of almost free
divisors in terms of differential forms.

For a reasonable description, the usual K\"{a}hler forms $\O^p_D$ (the 
definition is recalled at the start of Section \ref{def}) are not suitable:
there are too many of them. It is well known that if $x$ is a singular point
of a divisor $D$ then $\O^k_{D,x}$ has torsion for 
$\dim\ D + 1\geq k\geq \dim\ D -\dim D_{\mbox{sing}}$. 
The torsion gives important
geometrical information: indeed, if the $n$-dimensional
divisor $D$ has an {\it isolated} singularity
at $x$ then the torsion submodule of $\O^n_{D,x}$ has length equal to
the Tjurina number, and thus less than or equal to the Milnor number, 
the rank of the vanishing
homology (and equal to it if $D$ has a weighted homogeneous defining
equation). In this case one can think of torsion $n$-forms as forms which 
are idle
(evaluate to 0 on every n-vector) on $(D,x)$, but which come to life on 
the Milnor fibre of $(D,x)$. The torsion submodule of $\O^n_{D,x}$
thus encodes information on the vanishing 
geometry of the Milnor fibre. On a singular free divisor or almost 
free divisor $D$ the 
singular subspace
has codimension 1, and thus if $\dim\ D > 1$, all the modules
of K\"{a}hler forms (except for $\O^0$) will have torsion submodules of 
infinite length. In the deformations described by Damon, almost free divisors
become free, but are not smoothed, and thus infinite dimensional spaces of  
K\"{a}hler $k$-forms will wait in vain for their vanishing $k$-vectors. 

The obvious thing to do to remedy this is to kill the torsion forms on 
a free divisor $D$ (as free divisors are the stable objects in Damon's theory). 
Doing so, we obtain modules which we denote $\ch^k_D$. 
The definition of the 
modules $\ch^k_{D_0}$ of forms on an almost free divisor $D_0$ then 
follows by a standard 
universal procedure --- embed $D_0$ in a suitable free 
divisor $D$, as the fibre of an ambient submersion to a smooth space $S$, 
and then take the quotient of
the module of relative forms $\ch^k_{D/S}$ defined in the standard way from 
$\ch^k_D$.  
In fact Damon defines an almost free divisor $D_0\subset V$ 
as the pull back of a free divisor
$E\subset W$ by a map $i_0:V\to W$ which is transverse to $E$
outside $0$. 
In the spirit of Damon's definition, we can also define our complex
$\ch^\b_{D_0}$ as the quotient
of $\O^\b_V$ by the ideal generated by the pull-back via $i_0$ of a certain
ideal of $\O^\b_W$.
This is discussed in Section \ref{def}. 

This programme, which of course can be applied to any category
of spaces with a distinguished category of deformations, turns out to work
extremely well for free and almost-free divisors. In fact if $D$ is a free
divisor then 
$$\O^k_D/\mbox{torsion} = \O^k_X/h\O^k(\log D)$$
where $h$ is an equation for $D$ in $X$, and 
so, since all the modules $h\O^k(\log D)$ are free,
the $\ch^k_D$ are maximal Cohen-Macaulay $\OO_D$-modules. In Section
\ref{algebraic properties} we show that as a consequence, the 
$\ch^k_{D/S}$ enjoy
precisely the same depth properties as the modules $\O^k_{X/S}$ on deformations
of isolated complete intersection singularities (ICIS). In particular we
prove an analogue of the fact, recalled above, that for an isolated
hypersurface
singularity the torsion has length $\tau$:
\begin{theorem}\label{tor}
Let $f_0:(\CC^n,0)\to(\CC^p,0)$, with $n\geq p$ nice dimensions,
have finite ${\cal A}_e$-codimension, and let $D(f_0)$ be its discriminant.
Then the torsion submodule of $\ch^{p-1}_{D(f_0),0}$ has length equal to
the ${\cal A}_e$ codimension of $f_0$.
\end{theorem}
This is a consquence of \ref{torsion}(4) below.\\

\ni In Sections~\ref{using}, ~\ref{gm} and ~\ref{calc} below, we develop the
theory of the Gauss-Manin connection on the vanishing cohomology of
a deformation of an almost free divisor. The similarity with the classical
theory --- the Gauss-Manin connection on cohomology of the Milnor fibration
of an isolated hypersurface singularity
--- is evident; many of the same theorems hold, including coherence, 
and with almost the same proofs. The author does not claim much originality
in this latter part of the paper. However, what is remarkable is that the old
constructions
work so well in this new context.\\ 

\ni This paper has taken a rather long time to write, in the course of which
the author, and hopefully the paper, have profited from
many helpful conversations with
Francisco Castro and
Luis Narv\'{a}ez, to whom I am very grateful. Both the idea of studying the
Gauss-Manin connection on the discriminant of a family, and the definition
of $\ch^k_D$,  arose in conversation with them. I am also grateful to 
Jim Damon and Duco
van Straten for stimulating conversations on this topic, and to Mar\'{i}a
Aparecida Ruas, 
Farid Tari and Washington Luiz Marar for an invitation to lecture
on this topic in the 1998 S\~{a}o Carlos Singularities Workshop,
which led to what I believe are some useful improvements in the
presentation of this paper. I thank the referee for a very careful reading 
of the manuscript, and for a number of helpful suggestions.

Before the deadlines
passed me by, I optimistically intended the paper for the proceedings of the
1996 Oberwolfach conference in honour of E.Brieskorn on his 60th birthday;
notwithstanding the
time that has since passed, the paper remains respectfully dedicated to
him: imitation is the sincerest form of flattery

\section{Almost-free divisors and their deformations}
Let $E\subset W$ be a divisor in the manifold $W$, and let $V$ be a manifold.
The map $i:V\to W$ is {\it algebraically transverse} to $E$ at $w=i(v)$
if $d_vi(T_vV)+T^{\log}_wE=T_wW$. Here $T^{\log}_wE$ is the ``logarithmic'' 
tangent space to $E$ at $w$:
$$T^{\log}_wE=\{\chi(w):\chi\in\Der(\log E)_w\}.$$
For future reference we note that the notion of logarithmic tangent space
allows us to speak of logarithmic submersion and immersion, 
logarithmic critical value, etc. etc., as well as 
logarithmic ($=$ algebraic) transversality. 
By 2.12 of \cite{fischer} 
the logarithmic implicit function theorem holds (for maps with smooth
target). 
\begin{definition}(\cite{damonafd})  
The germ of hypersurface $D_0\subset V,0$ is an almost free 
divisor based on the
 germ of free divisor $E\subset W,0$ if there is a map $i_0:V,0\to W,0$ which is
algebraically transverse to $E$ except at $0$, such that $D_0=i_0^{-1}(E)$.
\end{definition}
If the map $i_0$ of the definition is algebraically transverse to $E$ then
$D_0$ is a free divisor; 
in this case it is not too hard to show that there is a
free divisor $E_0\subset W_0$ and integers $k$ and $\ell$ such that 
$(W,E)\simeq (W_0,E_0)\times\CC^k$ and 
$(V,D_0)\simeq (W_0,E_0)\times\CC^\ell$.

A 1-dimensional almost free divisor $D_0\subset V$ based on the free divisor
$E\subset W$ is also free in its own right, being 
a plane curve, even when the inducing map $i_0:V\to W$ is not algebraically 
transverse to $E$. Although this phenomenon is not inherently problematic, it
can be confusing, and therefore we now show
that it cannot occur for higher dimensional almost free divisors.

\begin{proposition}\label{gasp} Let $D_0\subset V$ be an almost free divisor 
based on the free divisor $E\subset W$, $D_0=i_0^{-1}(E)$. 
Then\\ 
(1)\ The sequence
\begin{equation}
0\to \Der(\log D_0)\ 
\hook\ \th_V\stackrel{ti_0}{\too} 
\frac{\th(i_0)}{i_0^\ast(\Der(\log E))}\to \frac{\th(i_0)}
{T{\cal K}_{E,e}i_0}\to 0
\label{seslog}
\end{equation}
(where $T{\cal K}_{E,e}i_0=ti_0(\theta_V)+i_0^\ast(\Der(\log E)$) is exact.\\
\ni (2)\ 
If $i_0$ is not algebraically 
transverse to $E$ then $\depth\ \Der(\log D_0)=2$ (we assume $\dim\ V\geq 2$).\\

\ni(3)\ (Jacobian criterion for freeness)
If $\dim\ V > 2$ then $D_0$ is a free divisor if and only if $i_0$ 
is algebraically transverse to $E$.
\end{proposition}
\Proof
First we prove that $i_0^\ast(\Der(\log E))$ is a free $\OO_V$-module.
Let $\xi_1,\cdots,\xi_n$ be a free basis for 
$\Der(\log E)$. We have an epimorphism 
$$\oplus^n_1 \OO_V\stackrel{\phi}{\too} i_0^\ast(\Der(\log E))$$
with $\phi(e_j)=\xi_j\circ i_0$. Suppose 
$(a_1,\cdots,a_n)\in \ker\ \phi$. Then $\sum_j a_j(\xi_j\circ i_0)=0$.
However, if $i_0(x)\notin E$ then the vectors 
$\xi_1(i_0(x)),\cdots,\xi_n(i_0(x))$
are linearly independent, and so all $a_j(x)$ vanish. That is, $\supp(a_j)
\subset D_0$. As $D_0$ is a divisor, all the $a_j$ must be identically $0$,
and $i_0^\ast(\Der(\log E))$ is free.\\

\ni (1) and (2) Consider the short exact sequence
$$0\to i_0^\ast(\Der(\log E))\to\th(i_0)\to 
\th(i_0)/i_0^\ast(\Der(\log E))\to 0.$$
The right-most module is supported only
on $D_0$ and thus is not free over $\OO_V$, so its depth is less than 
the depth of $\th(i_0)$, which is of course equal to $\dim\ V$ as $\th(i_0)$ is
a free $\OO_V$-module. As $i_0^\ast(\Der(\log E))$ is also free, the depth lemma
shows that $\depth\ \th(i_0)/i_0^\ast(\Der(\log E))= \dim\ V-1$. 

The exact sequence 
$$0\to K\to \th_V\stackrel{ti_0}{\too} \th(i_0)/i_0^\ast(\Der(\log E))\to 
\th(i_0)/T{\cal K}_{E,e}i_0\to 0$$
(where $K$ and $T{\cal K}_{E,e}i_0$ are defined by the sequence) 
breaks into the two short exact 
sequences
$$0\to \th_V/K\to \th(i_0)/i_0^\ast(\Der(\log E))\to \th(i_0)/
T{\cal K}_{E,e}i_0\to 0$$
and
$$0\to K\to\th_V\to \th_V/K\to 0.$$
If $\th(i_0)/T{\cal K}_{E,e}i_0\neq 0$, then its depth is $0$ (for it is
supported only at $0$). As\\ 
$\depth\ \th(i_0)/i_0^\ast(\Der(\log E))= \dim\ V-1>0$,
the depth lemma applied to the
first of these two short exact sequences 
implies that $\depth\ \th_V/K=1$. The second short exact sequence now implies
that $\depth\ K=2$.

A similar argument shows that
if $\th(i_0)/T{\cal K}_{E,e}i_0=0$ then $K$ must be free.

Now $K\subset \Der(\log D_0)$, for if $\chi\in K$
then for all $x\in D_0$, $di_{0\ x}(\chi(x))\in T^{\log}_xE$; if $x$ is a
regular point of $D_0$ this means that $i_0(x)$ is a 
regular point of $E$ and that
$\chi(x)\in (d_xi_0)^{-1}(T_{i_0(x)}E)=T_xD_0$, and thus 
$\chi\in\Der(\log D_0)$. Since
$i_0$ is algebraically transverse to $E$ outside $0$, $K$ and $\Der(\log D_0)$
coincide outside $0$. Thus we have a short exact sequence
$$0\to K\to \Der(\log D_0)\to \Der(\log D_0)/K\to 0$$
in which the right-most module, if not $0$, has dimension $0$ and thus 
depth $0$. This is impossible,
again by the depth lemma,
since $\depth\ K\geq 2$. Thus $K=\Der(\log D_0)$.\\

\ni (3) ``if'' is clear; here we prove the converse. 
If $\th(i_0)/T{\cal K}_{E,e}i_0\neq 0$ then by 
(3) $\depth\ \Der(\log D_0) =
2 < \dim\ V$, and thus $D_0$ cannot be a free divisor.  \eop 
\begin{example} {\em The divisor $D_0$ shown in the Introduction is
an almost free divisor. It is the preimage of the free divisor $E$
consisting of the union of the four co-ordinate hyperplanes in $\CC^4$,
under the inducing map $i_0(x_1,x_2,x_3)=(x_1,x_2,x_3,x_1+x_2+x_3)$. Since
$i_0$ is not algebraically transverse to $E$ at $0$ (for $T^{\log}_0E=\{0\}$),
it follows from \ref{gasp}(3) that $D_0$ is not free at $0$. Since $D_0$ is 
free outside $0$, it also follows that $i_0$ is algebraically transverse to
$E$ outside $0$ (though of course this is obvious anyway), so $D_0$ is an
almost free divisor.}
\end{example}

The deformations of almost free divisors that
we consider arise as follows: 
beginning with an almost free divisor $D_0\subset V$ obtained by 
pulling-back the free divisor $E\subset W$ by the map $i_0:V\to W$, 
we consider deformations $i:V\x S\to W$ of
$i_0$, and fibre $D:=i^{-1}(E)$ over the base $S$ of the deformation.
When
necessary we will refer to deformations of this type as {\it admissible}.
Unless otherwise specified, {\it all deformations of almost free divisors 
that we consider from now on will be admissible}. 
A deformation $D\to S$ is {\em free} if $D=i^{-1}(E)$ with $i$ 
algebraically transverse to $E$, is {\it versal} if $i$ is a ${\cal K}_E$-versal
deformation of $i_0$ (see 2.6 below),
and {\it frees $D_0$} if for generic $u\in S$
the map $i_u=i(\_,u)$ is algebraically transverse to $E$. The set 
of points $u\in S$ for which $i_u$ is {\it not} algebraically transverse
to $E$ is the {\it logarithmic discriminant} of the deformation, and is denoted
${\cal B}$ (since in many cases it is the bifurcation set of a deformation).
When the deformation is miniversal, the logarithmic discriminant is 
sometimes called the {\it ${\cal K}_E$-discriminant} of $D_0$ (or of
$i_0$).
\begin{theorem/definition}
\cite{dm}\cite{damonafd} Let $D_0$ be a $p$-dimensional almost free divisor
based on the free divisor $E$. If $\pi:D\to S$ is a deformation of
$D_0$ which frees $D_0$, then for 
$u\in S\setminus {\cal B}$, 
$D_s=i_u^{-1}(E)$ has the homotopy type of a wedge of $p$-spheres. The number
of these is independent of the choice of deformation; it is called 
the {\em singular Milnor number} of $D_0$, and denoted $\mu_E(D_0)$. The
space $D_s$ is a {\em singular Milnor fibre} of $D_0$; up to homeomorphism
it is independent of the choice of deformation. 
\end{theorem/definition}
If the logarithmic
version of Sard's theorem held for free divisors,
then every free deformation of $D_0$ would free it. 
However, it does not hold:

\begin{example} {\em (F. J. Calder\'{o}n, \cite{tesis})
Let $D\subseteq\CC^3$ be the set defined by the equation
$xy(x-y)(x+\lambda y)=0$, and let $\pi(x,y,\lambda)=\lambda$.
One checks that $D$ is a free
divisor.
The fibre of $\pi$ over $\lambda$ is
a union of four lines through $(0,0)$ in $\CC^2$. As $\lambda$ varies,
the cross-ratio of these four lines varies, and thus the family is
not trivial. In particular, there is no logarithmic vector field with
non-vanishing $\p/\p \lambda$-component at points $(0,0,\lambda)$, and so
every point in $\CC$ is a logarithmic critical value of $\pi$.}
\end{example}
In \cite{damonagain} Jim Damon describes a large class of related examples; 
in particular the total space of the Hessian deformation of any non-simple
weighted homogeneous plane curve singularity is a free divisor. 
Calder\'{o}n's example falls in this class, of course.\\

\ni It is probably sensible at this point to describe in more detail
the family of examples we
are most interested in. These are the discriminants of map-germs
of isolated instability (i.e. finite ${\cal A}_e$ codimension) 
$f_0:U=\CC^n,0\to\CC^p,0=V$ (with $n\geq p$), and those of 
their deformations which
arise from
deformations of the germ $f_0$.
Any map $f_0:U\to V$
can be obtained by
transverse pull-back $i_0$ from a stable map $F$,
as in the diagram
$$\begin{array}{rcl}
{\cal U}&\stackrel{F}{\too}&{\cal V}\\
j_0\uparrow&&\uparrow\ i_0\\
U&\stackrel{f}{\too}&V
\end{array}
$$
In this case the discriminant $D(f_0)$ of $f_0$ is the pull-back by $i_0$ of the
discriminant $D(F)$ of $F$. 
In \cite{damonKV}, 
Damon shows that stability of $f$ is equivalent to the algebraic
transversality of $i_0$ to $D(F)$. We reprove this in Section \ref{maps} below.
Thus, if $f_0$ has isolated instability then
$D(f_0)$ is an almost free divisor based on $D(F)$. 

Any unfolding $G:U\x S\to V\x S,\ g(x,s)=(g_s(x),s)$\ \  of
$f_0$ can also be obtained from $F$ by transverse pull-back, by a deformation
$i$ of $i_0$, and thus the fibration $D(G)\to S$ is an admissible deformation
of $D(f_0)$. Its fibre over $s\in S$ is the discriminant $D(g_s)$ of $g_s$.
Now suppose that $G$  is a stable unfolding of $f_0$, 
and that the logarithmic discriminant 
$\cal B$ 
of the projection $\pi:D(G)\to S$
is a proper subset of $S$ (so that $G$ is a stabilisation of $f_0$).
For $u\notin {\cal B}$, $f_u$ is a stable perturbation
of $f_0$. The singular Milnor number $\mu_{D(F)}(D(f_0))$ of 
$D(f_0)$ is the rank of
the vanishing homology of the discriminant of
a stable perturbation of $f_0$; in \cite{dm} it is denoted by
$\mu_\Delta(f_0)$ and called the discriminant Milnor number of $f_0$.
The main theorem 
of \cite{dm} is that provided $(\dim U,\dim V)$ are in Mather's
"nice dimensions", 
$\mu_{\Delta}(f_0)$) satisfies 
$$\mu_\Delta(f_0)\geq {\cal A}_e\mbox{codimension of}\ f_0$$
with equality in case $f_0$ is weighted homogeneous. We give a fuller account
of the link
between the theory we construct here and the theory of right-left equivalence
of map-germs, in Section~\ref{maps} below. For now, we hope 
it will serve as a motivating example. \\

Now recall from \cite{damondef} the notion of ${\cal K}_E$
equivalence for almost free divisors based on the free divisor $E$ - or, more
precisely, ${\cal K}_E$-equivalence of the maps inducing them. 

\begin{definition} 1. Let $V$ and $W$ be smooth, and let $E\subset W$ be a 
closed complex subspace. 
Then ${\cal K}_E$ is
the subgroup of $\cal K$ consisting of pairs $(\Phi,\phi)\in 
\mbox{Diff}(V\x W)\x \mbox{Diff}(V)$ such that
\begin{itemize} 
\item
$p_1(\Phi(v,w))=\phi(v)$ for $(v,w)\in V\x W$, 
where $p_1:V\x W\to V$ is projection;
\item
$\Phi(v,0)=(\phi(v),0)$ for $v\in V$;
\item
$\Phi(V\x E)=V\x E$.
\end{itemize}
2. Two germs $i_0,j_0:V\to W$ are ${\cal K}_E$-equivalent if there exists
$(\Phi,\phi)\in {\cal K}_E$ such that $\Phi\circ \mbox{gr}_{i_0}
=\mbox{gr}_{j_0}\circ
\phi$, where $\mbox{gr}_{i_0}$ and $\mbox{gr}_{j_0}$ 
are the graph embeddings of
$i_0$ and $j_0$, $v\mapsto (v,i_0(v))$ and $v\mapsto (v,j_0(v))$ respectively.\\
\end{definition}

\ni Evidently if $i_0$ and $j_0$ are ${\cal K}_E$-equivalent then 
$i_0^{-1}(E)$ and
$j_0^{-1}(E)$ are isomorphic varieties.\\

\ni The ${\cal K}_{E,e}$ tangent space $T{\cal K}_{E,e}i_0$ 
is equal to
$$ti_0(\t_V)+i_0^\ast\Der(log E);$$ 
the ${\cal K}_{E,e}$ normal space is then
$$\frac{\t(i_0)}{ti_0(\t_V)+i_0^\ast\Der(\log E)}.$$
By Nakayama's lemma, it is supported precisely on the points where $i_0$
fails to be algebraically transverse to $E$.
The deformation $i$ of $i_0$ is ${\cal K}_{E,e}$ versal if every other
deformation of $i_0$ is parametrised ${\cal K}_E$ isomorphic, 
in an appropriate sense, to one induced from $i$. The sub-index ``$e$''
here indicates that in the parametrised equivalence, the parametrised family
$\phi_s$ of diffeomorphisms of the domain $V$ of $i+-$ are not
required to fix $0\in V$ when $s\neq 0$. Recquiring that $\phi_s(0)=0$
for all $s\in S$ gives rise to a stricter notion of versality, namely
${\cal K}_E$-versality. Note that our terminology differs slightly from
Damon's at this point: he does not consider the stricter notion, and
uses the term ${\cal K}_E$-versal where we use ${\cal K}_{E,e}$ versal.
See Section 3 of \cite{wallsurvey} for a discussion of this point.

Damon shows in \cite{damondef} that the usual infinitesimal criterion for
versality holds: $i$ is ${\cal K}_{E,e}$-versal if and only if 
$$T{\cal K}_{E,e}i_0+\CC\langle\dot{i}_1,\cdots,\dot{i}_d\rangle = \t(i_0)$$
where $\dot{i}_j$ is the restriction to $V\x \{0\}$ of $\p i/\p s_j$, 
the $s_j, j=1,\cdots,d$ being coordinates on $S$. 

From the existence of ${\cal K}_{E,e}$-versal deformations it follows that
if $\pi:D_1\to S_1$ and $\pi_2:D_2\to S_2$ are both admissible deformations
of the almost free divisor $D_0$ based on the free divisor $E$, and both 
free $D_0$, then the generic fibres of $\pi_1$ and $\pi_2$ are
homeomorphic. For up to isomorphism 
both deformations are induced by base change from a versal
deformation, and as the base space of the latter is smooth, and the
discriminant is a proper analytic subvariety, it follows that any
two generic fibres are homeomorphic.

\begin{lemma}\label{iso} Let $D_0\subset V$ be an almost 
free divisor based on the 
free divisor $E\subset W$, with $i_0:V\to W$ the inducing map. 
Let $i:V\times S$ be a deformation of
$i_0$, 
algebraically transverse to $E$, let $D=i^{-1}(E)$, and
let $\pi:D\to S$ be the resulting deformation of 
$D_0$. Then 
$$N{\cal K}_{E,e}i_0 \simeq \frac{\t(\pi)}{t\pi(\Der(\log D))+\pi^\ast(m_{S,0})
\t(\pi)}.$$
\end{lemma}
\Proof We have
$$\frac
{\t(\pi)}{t\pi(\Der(\log D))+\pi^\ast(m_{S,0})
\t(\pi)}
\simeq
\frac
{\theta_{V\times S}}
{\theta_{V\times S/S}+\Der(\log D)+m_{S,0}\theta_{V\times S}}$$
$$\simeq
\frac
{\theta(j)}
{tj(\theta_V)+j^\ast(\Der(\log D))}\ =\ N{\cal K}_{D,e}j,$$
where $j:V\to V\times S$ is the inclusion $v\mapsto (v,0)$. Now
$i_0=i\circ j$, and because $i$ is algebraically transverse to $E$,
$ti:\theta(j)\to\theta(i_0)$ induces an isomorphism $N{\cal K}_{D,e}j
\simeq N{\cal K}_{E,e}i$. This is proved in Proposition 1.5 of
\cite{damonKV}, but for the reader's convenience we now sketch the argument.
It is easiest to see if we assume that 
$T^{\log}_0E=0$; for then $i$ must be a submersion, and choosing coordinates
with respect to which it is a projection, the isomorphism is a straightforward
calculation. 

In the general case, if $T^{\log}_0E$ is $k$-dimensional, then the pair
$(W,E)$ is isomorphic in some neighbourhood to a product $(W_0,E_0)\times
\CC^k$. Let $p:W\to W_0$ be projection, using this product structure. Then
the calculation just mentioned shows that
$N{\cal K}_{E,e}i_0\simeq N{\cal K}_{E_0,e}p\circ i_0$, and also that
the submersion $p\circ i$ induces an isomorphism $N{\cal K}_{D,e}j\simeq
N{\cal K}_{E_0,e}p\circ i_0$.
\eop
In the light of this lemma, we refer to the quotient $\theta(\pi)/t\pi(\Der(\log D))$ as $T^{1,\log}_{D/S}$, and to $\theta(\pi)/t\pi(\Der(\log
 D))+m_{S,0}\theta(\pi)$ as $T^{1,\log}_{D_0}$. The reader will recognise the
similarity to the definition of $T^1_{X/S}$ and $T^1_{X_0}$ for a deformation
$f:X\to S$ of an ICIS $X_0$, cf \cite{looijenga} Chapter 6.

We extend the definition of $T^{1,\log}_{D/S}$
to the case where $\pi:D\to S$ is
not necessarily free:

\begin{definition} Let $\pi:D\to S$ be an admissible deformation, and let
$\rho:{\cal D}\to S\x T$ be a free extension of $\pi$. Then 
$$T^{1, \log}_{D/S}=\frac{\t(\rho)}{t\rho(\Der(\log {\cal D}))+(t)\t(\rho)},$$
where $(t)$ is the ideal defining $S\x \{0\}$ in $S\x T$.
\end{definition}
It is easily seen that this definition is independent of the choice of 
free extension $\rho$, by means of the following lemma, which is
proved by standard singularity theory methods: 
\begin{lemma}\label{welldef}
 Let $E\subset W$ be a divisor, let $i_0:V,0\to W,0$, and
suppose that $i:V\times T,(0,0)\to W,0$ and
$j:V\times T,(0,0)\to W,0$ are deformations of $i_0$, both of them 
logarithmically transverse to $E$. Then $i$ and $j$ are ${\cal K}_E$-equivalent,
as maps, by an equivalence which restricts to the identity on $V$.\eop
\end{lemma}
The analogy with isolated complete intersection singularities continues:
\begin{proposition}\label{cm} Let $\pi:D\to S$ be an admissible
deformation of the almost free divisor $D_0$, and 
suppose that $\pi$
frees $D_0$. Then
$T^{1,\log}_{D/S}$ is a Cohen-Macaulay module of
dimension $\dim S -1$.
\end{proposition}
\Proof If $\pi$ is a free deformation, then the result
 follows from a theorem of Buchsbaum and Rim (\cite{br}
Corollary 2.7); for we have an exact sequence
$$\Der(\log D)\to\t(\tilde{\pi})\to T^{1,\log}_{D/S}\to 0$$
(where $\tilde{\pi}:V\times S\to S$ is the obvious extension of $\pi$)
in which $\Der(\log D(F))$ is free of rank $\dim D+1$ and $\t(\pi)$ is
free of rank $\dim S$. The theorem cited
implies that $\dim T^{1,\log}_{D/S}\geq \dim S - 1$, and that if this dimension
is attained then $T^{1,\log}_{D/S}$ is Cohen-Macaulay. Now the hypothesis that
$D_0$ be almost free implies that $T^{1,\log}_{D/S}$ is finite over $S$,
while the hypothesis that $\pi$ frees $D_0$ implies that
$\dim(\supp(\pi_\ast(T^{1,\log}_{D/S}))<\dim S$.  Hence
$\dim T^{1,\log}_{D/S}= \dim S - 1$, and (1) follows.\\ 

\ni The general case follows from this special case by a standard
dimensional argument.
Let $\rho:{\cal D}\to S\x T$ be a free extension of $\pi$, with $\dim\ T=d$.
Then  because
$$T^{1,\log} _{ D/S}=(T^{1,\log}_{
{\cal D}/S\x T})/(t_1,\cdots,t_d)T^{1,\log} _{{\cal D}/S\x T}$$ and 
$\dim\ T^{1,\log} _{ D/S}=\dim\ (T^{1,\log}_{
{\cal D}/S\x T})-d$, $T^{1,\log} _{ D/S}$ also is Cohen-Macaulay.\\
\eop

\bc 
\section{Differential forms on free and almost free divisors}\label{def}
\ec

From now on, when we speak of a free divisor the letter $h$ will always denote
its reduced equation.

In this section we define and study the improved version of K\"{a}hler forms
on free divisors mentioned in the introduction.
First we recall (e.g. from \cite{looijenga}) the standard defintion of the
$\OO_X$-module $\O^k_X$ of K\"{a}hler $k$-forms on a subvariety $X$ of a 
smooth space $V$. Suppose that the ideal of functions vanishing on $X$
is generated by $f_1,\ldots,f_n$; then $\O^k_X$ is defined to be
$$\frac{\O^k_V}{(f_1,\ldots,f_n)\O^k_V+\sum_idf_i\wedge\O^{k-1}_V}.$$
The exterior derivative $d:\O^k_V\to\O^{k+1}_V$ passes to the quotient to
give an exterior derivative $\O^k_X\to\O^{k+1}_X$; the complex $\O^\b_X$,
equipped with this exterior derivative, is naturally isomorphic to the
complex of holomorphic differential forms on $X$ defined using charts, 
when $X$ is smooth.

Now let $D\subset V$ be a free divisor. 
Recall  K.Saito's definition  in \cite{saito} of the sheaf 
$\Omega^k(\log D)$ 
as the 
$\OO_V$-module of meromorphic differential forms $\o$ on $V$ such that
$h\o$ and $h\ d\o$ are regular.  
Clearly $\Omega^\b_V(\log D)$ is a complex;
note that $h\O^\b_V(\log D)$ is a subcomplex of $\O^\b_V$, so that
the quotient $\O^\b_V/h\O^\b_V(\log D)$ is also a complex. This complex is
our replacement for the complex of K\"{a}hler forms on $D$.
\begin{definition} Let $D$ be a free divisor. We define 
$$\ch^k_D=\frac{\O^k_W}{h\O^k(\log D)}.$$
\end{definition}
Since 
$(dh)/h\wedge \Omega^{k -1}_V+\Omega^k_V
\subseteq \O^k_V(\log D)$, $\ch^k_D$
 is a quotient of $\O^k_D$.
\vskip\baselineskip
\begin{example} {\em Let $D$ be the normal crossing divisor $z_1\ldots z_p=0$.
Then $\O^1_V(\log D)$ is the free $\OO_V$-exterior algebra on generators 
$dz_i/z_i, i=1,\cdots,p$, and so 
$h\O^k_V(\log D)$ is generated over $\OO_V$ by the forms $z_Jdz_I$, where
$|I|=k$, 
$J\cup I=\{1,\cdots,p\}$ and $J\cap I=\emptyset$. It is easy to see that
if $n:\bar{D}\to D$ is the normalisation of $D$, and $i:D\to \CC^{n+1}$ is 
inclusion, then  
$h\O^k_V(\log D)=\ker(i\circ n)^\ast:
\O^k_V\to\O^k_{\bar{D}}$, and so $\ch^k_D 
\hookrightarrow n_\ast
\O^k_{\bar{D}}$. By contrast, $\O^k_D$ is larger, 
and the kernel
of $n^\ast:\O^k_D\to\O^k_{\bar D}$ is non-trivial:
for example,
the forms $(h/x_i)dx_i$ are not equal to $0$ in $\O^1_D$.} 
\end{example}
\begin{lemma}\label{absbits}
Let $D$ be a divisor in the $p+1$-dimensional complex 
manifold $V$, with local defining equation $h$. Then\\

\ni (1) $\ch^0_D=\OO_D$.\\

\ni (2) $\ch^{p+1}_D=0$.\\

\ni (3) For $k$ with $0\leq k\leq p$, 
the support of $\ch^k_D$ is equal to $D$.\\

\ni (4) If $x$ is a smooth point of $D$, then $\ch^k_{D,x}=\O^k_{D,x}$.\\ 

\ni (5) If $D$ is a free divisor, then 
the depth (as module over $\OO_V$ and $\OO_D$) of 
$\ch^k_D$ is equal to $\dim D$, for $0\leq k\leq p$,
(so that
$\ch^k_D$ is a maximal Cohen Macaulay $\OO_D$-module).\\

\ni (6) If the germ of $D$ at $x$ is quasihomogeneous, then the complex
$$0\to\CC\to\OO_D=\ch^0_D\to\ch^1_D\to\cdots\to\ch^p_D\to 0$$
is exact at $x$. 
\end{lemma}
\Proof
(1) and (3) are obvious.\\

\ni (2) follows from the fact that $\O^{p+1}(\log D)=(1/h)\O^{p+1}_V$.\\

\ni (4). If $y$ is a smooth point of $D$, then 
$\O^k(\log D)_y=(dh/h)\wedge \Omega^{k-1}_{V,y}+\Omega^k_{V,y}$.
Hence $h\Omega^k(\log D)=dh\wedge\Omega^{k-1}_{V,y}+h\O^k_{V,y}$,
and $\O^k_{V,y}/h\O^k(\log D) = \O^k_{D,y}$.\\

\ni (5). $\O^k_V/h\O^k_V(\log D)$ has an $\OO_V$-free resolution 
$0\to h\O^k(\log D)\to \O^k_V\to\O^k_V/h\O^k(\log D) \to 0$, so its
projective dimension is less than or equal to $1$. Since its support
is $D$, it must have projective dimension $1$ and thus depth $p-1$.\\

\ni (6). The proof here is a variant of the elementary proof
that the usual complex 
$$0\to\CC\to\OO_D\to\O^1_D\to\cdots\to\O^{p-1}_D\to \O^p_D\to 0$$ 
is exact, which can be found, for example, as Lemma 9.9 in 
\cite{looijenga}. 
We use the local $\CC^\ast$-action centred at $y$. Let $\chi_e$ be a local
Euler field, vanishing at $y$. Then for any homogeneous form $\omega\in
\O^k_{V,y}$ of weight $\ell$, we have $L_{\chi_e}(\omega)=\ell \o$.
Note Cartan's identity 
$L_{\chi_e}(\o)=\i_{\chi_e}(d\o)+d\i_{\chi_e}(\o)$. Now suppose that
$d\o\in h\O^k(\log D)$, and for brevity write $d\o=h\sigma$. Then
$\ell\o=h\i_{\chi_e}(\sigma)+d\i_{\chi_e}(\o)
\ \in h\O^k(\log D)_y+d\O^{k-1}_{V,y}$.  Thus, unless $\ell=0$, the class
of $\o$ in $\O^k_{V,y}/h\O^k(\log D)_y$ is a coboundary. It follows
that if we define $H:\O^\b_{V,y}/\O^\b(\log D)\to \O^{\b-1}_{V,y}/
h\O^{\b-1}(\log D)_y$ to be $(1/\ell)\i_{\chi_e}$ on the weight $\ell$
subspace (and extend continuously to the completion), then $H$ gives
a homotopy equivalence between $\O^\b_{V,y}/\O^\b(\log D)_y$ and its weight
zero subspace. This weight 0 subspace is just $\CC$ in degree 0, since 
our $\CC^\ast$ action is good (i.e. the weights of all the variables are
positive).
\eop
\begin{corollary}
Let $D\subset V$ be a free divisor. For each $k$, with $0\leq k\leq \dim V$,
the torsion submodule of $\O^k_D$
is equal to $h\O^k(\log D)/h\O^k_V+dh\wedge \O^{k-1}_V$, and thus
$\ch^k_D=\O^k_D/\mbox{torsion}$.
\end{corollary}
\Proof
As $\ch^k_D$ is Cohen-Macaulay it is torsion free, and thus
$$T\O^k_D\subseteq h\O^k(\log D)/h\O^k_V+dh\wedge\O^{k-1}_V.$$ However
the right hand side here is torsion, as it is supported only on
$D_{\mbox{sing}}$ (cf. \ref{absbits}.5.)
\eop
The fact that on a free divisor $\ch^\b_D$ 
is a complex of maximal
Cohen-Macaulay modules on $D$ is crucial to its applications. 
In this respect it improves on the usual complex $\O^\b_D$; in fact 
\begin{proposition}\label{bloom} If $D$ is a free divisor and singular at $x$
then for $1\leq k\leq \dim\ D$, 
$\depth\ \O^k_D = \dim\ D-1$.
\end{proposition}
\Proof This follows by applying the depth lemma to each of the 
three short exact sequences 
$$0\to \O^k_V\to\frac{\O^k(\log D)}{\O^k_V}\to 0,$$
$$0\to \frac{\O^{k-1}_V}{h\O^{k-1}(\log D)}\stackrel{(dh/h)\wedge}
{\longrightarrow}\ \frac{\O^k(\log D)}{\O^k}\to
\frac{\O^k(\log D)}{(dh/h)\wedge \O^{k-1}+\O^k}\to 0$$
and
$$0\to\frac{\O^k(\log D)}{(dh/h)\wedge \O^{k-1}+\O^k}
\stackrel{\times h}{\longrightarrow}
\frac{\O^k}{dh\wedge\O^{k-1}+h\O^k}\to\frac{\O^k}{h\O^k(\log D)}\to 0$$
together with~\ref{con} below. For the first shows that $\O^k(\log D)/\O^k_V$
has depth $p:=\dim\ D$, and so the second shows that its right-hand member has depth $p-1$ (since by \ref{absbits}(4) is is supported only on $D_{\mbox{sing}}$.
and by \ref{con} below it is not $0$). The third sequence has $\O^k_D$ as its
middle term; the depth lemma says that its depth is at least as great as the
minumum of the depths of the two outer terms, and thus at least $p-1$. It
also says that if it is greater than this minumum then the first term has
depth one more than the last. Since the first term has depth less than the last,
this cannot happen.
\eop
\begin{lemma}\label{con} If $x$ is a singular point on the reduced free divisor
$D\subseteq \CC^{p+1}$ then for $0 < k \leq p+1$,
$dh\wedge\O^{k-1}_{\CC^{n+1},x}+h\O^k_{\CC^{n+1},x}$ 
is strictly contained in $h\O^k(\log D)_x$.
\end{lemma}
\Proof 
Suppose $\O^1(\log D)_x=\O^1_{\CC^{p+1},x}+(dh/h)\wedge\O^0_{\CC^{p+1},x}$. 
Then a free basis for
$\O^1(\log D)_x$ can be extracted from the list $dx_1,\cdots,dx_{p+1},dh/h$.
If $\o_1,\cdots,\o_{p+1}$ is a free basis,
then $\o_1\wedge\cdots\wedge\o_{p+1}=dx_1\wedge\cdots\wedge dx_{p+1}/h$,
up to multiplication by a unit
(\cite{saito}). Thus, $dh/h$ must be a member of any free basis extracted
from this list, and so without loss of generality, we can suppose that
$dx_1,\cdots,dx_p,dh/h$ is a free basis. Since $dx_{p+1}$ is then an 
$\OO_{\CC^{p+1}}$-linear combination of $dx_1,\cdots,dx_p,dh/h$, a 
calculation shows
that $\p h/\p x_{p+1}$
is a unit, and thus that $D$ is non-singular at $x$.

Let $\o_1,\ldots,\o_{p+1}$ be a free basis for $\O^1(\log D)$. 
Then $\O^\b(\log D)_x$ is the free $\OO_{\CC^{p+1},x}$
exterior algebra over
$\OO_{\CC^{p+1},x}$ generated by $\o_1,\ldots,\o_{p+1}$, while the algebra
$\O^\b_{\CC^{p+1},x}+(dh/h)\wedge\O^{\b -1}_{\CC^{p+1},x}$ is a quotient
(because there are relations in general) of the free $\OO_{\CC^{p+1},x}$
exterior algebra generated by $dx_1,\ldots,dx_{p+1},dh/h$.
Thus, if $\O^1_{\CC^{p+1},x}+(dh/h)\wedge\O^0_{\CC^{p+1},x}$ is strictly
contained in $\O^1(\log D)_x$ then for each $k$ with $0 < k\leq p+1$,
$\O^k_{\CC^{p+1},x}+(dh/h)\wedge\O^{k-1}_{\CC^{p+1},x}$ is strictly
contained in $\O^k(\log D)_x$.
\eop
Recall that $\theta_D=\mbox{Hom}_{\OO_D}(\Omega^1_D,\OO_D)$ is equal to the
restriction to $D$ of $\Der(\log D)$. The next proposition characterises
$\ch^1_D$ as the double dual of $\O^1_D$.
\begin{proposition}\label{pairing}
Let $D\subset V$ be a free divisor; then evaluation of forms on vector
fields induces a perfect pairing $\ch^1_D\times\theta_D\too \OO_D$.
\end{proposition}
\Proof The result is well known when $D$ is non-singular, so we assume $D$
is singular.
Define $N$ by the short exact sequence
$0\to\Der(\log D)\to\theta_V\to N\to 0$. It is a maximal Cohen-Macaulay
$\OO_D$ module. Because $\Der(\log D)\supset h\theta_V$,
restricting to $D$ gives the exact sequence
$0\to\theta_D\to\theta_V|_D\to N\to 0$,
and now dualising with respect to $\OO_D$ gives an exact sequence
$$0\to N^\vee\to (\theta_V|_D)^\vee 
\to \theta_D^\vee\to \mbox{Ext}^1_D(N,\OO_D)
\to 0.$$
Because $N$ is a maximal Cohen-Macaulay module on a singular hypersurface,
it has a 2-periodic $\OO_D$-free resolution (see section 6 of \cite{eis}), 
and this
remains exact on dualising with respect to 
$\OO_D$; hence $\mbox{Ext}^1_D(N,\OO_D)=0$. Thus
$\O^1_V|_D$ maps onto $\theta_D^\vee$, and hence so does $\O^1_V$.
The kernel of the epimorphism $\Omega^1_V\to\theta_D^\vee$ is 
$\{\o\in\O^1_V:\o(\chi)\in(h)\ \forall \chi\in\Der(\log D)\}$. 
Denote this by $K$. Since 
evaluation of forms on vector fields gives an isomorphism $\O^1(\log D)
\simeq \mbox{Hom}_V(\Der(\log D),\OO_V)$, $K$ is equal to $h\O^1(\log D)$.
Thus, $(\theta_D)^\vee$ is equal to $\ch^1_D$. 
Finally, since $\ch^1_D=(\O^1_D)^{\vee\vee}$, its dual must be $(\O^1_D)^\vee$,
i.e. $\theta_D$.
\eop
We now define differential forms on almost free divisors. Rather than
beginning with the procedure outlined in the introduction, we give an
alternative definition, more in the spirit of Damons's definition of
almost free divisor, in terms of the inducing map $i_0$.
The two approaches are shown to be equivalent in \ref{samediff} below.
\begin{definition}
Let $E\subset W$ be a free divisor,  let $i_0:V\to W$ be a map and 
let $i:V\times S\to W$ a deformation of $i_0$. Let
$D_0=i_0^{-1}(E)$ and $D=i^{-1}(E)$. 
We set
\be
\item
$$\ch^\b_{D_0}=\frac{\O^\b_V}{\langle i_0^\ast(h\O^\b(\log E))\rangle}$$
\item
$$\ch^\b_{D/S} = \frac{\O^\b_{V\times S}}{\langle i^\ast(h\O^\b(\log E))\rangle
+\sum ds_i\wedge \O^{\b-1}_{V\times S}}.$$
\ee
\end{definition}
Here $\langle i_0^\ast(h\O^\b(\log E))\rangle$ and 
$\langle i^\ast(h\O^\b(\log E))\rangle$ are the ideals in the 
exterior algebras $\O^\b_V$ and $\O^\b_{V\times S}$ generated by the forms 
$i_0^\ast(\o)$ and $i^\ast(\o)$, for
$\o\in h\O^\b(\log E)$. They are of course subcomplexes.
\begin{proposition}\label{samediff} 
Suppose that $i:V\times S\to W$ is a deformation
of $i_0:V\to W$. 
Let $h_E$ be a reduced equation for $E$, so that $h_D=h\circ i$ is a
reduced equation for $D=i^{-1}(E)$ .
Then
\be
\item
$$\ch^k_{D_0}=
\frac{\ch^k_{D/S}}
{\sum s_i\ch^k_{D\times S}}$$
(so the construction of $\ch^\b$ commutes with restriction).
\item
If $I$ is algebraically transverse to $E$ then
$$\ch^k_{D/S}=\frac{\O^k_{V\times S}}{h_D\O^k(\log D)
+\sum_ids_i\wedge\O^{k-1}_{V\times S}}$$
\ee
\end{proposition}
\Proof The first statement is obvious. For the second, we need a lemma.
\begin{lemma} Let $p:U\to W$ be a submersion, let $E\subset W$ be a 
divisor and let $D=p^{-1}(E)$. Then 
$$\O^\b(\log D)=\langle p^\ast(\O^\b(\log E))\rangle.$$
\end{lemma}
\Proof 
Inclusion of the right hand side in the left follows from the obvious fact that
$p^\ast(\O^\b(\log E))\subset \O^\b(\log D)$. 
To prove the opposite inclusion
note that after a diffeomorphism we may assume
$U=W\times T$ and $D=E\times T$, and that $p:U\to W$ is just the projection
$W\times T\to W$. Let $q:W\to W\times T$ be the inclusion $w\mapsto (w,0)$,
and let $t_1,\cdots,t_n$ be coordinates on $T$; then if
$\o\in\O^\ell(\log D)$, we have 
$$
\begin{array}{lll}
\o(x,t)-\o(x,0) 
&
= 
&\int_0^1 \frac{d}{ds}\o(x,st)ds
\\
&
=
&
\int_0^1\sum_i t_i\frac{\partial}{\partial t_i}\o(x,st)ds
\\
&
=
&
\sum_it_i\int_0^1\frac{\partial}{\partial t_i}\o(x,st)ds.
\end{array}
$$
Now each form $\int^1_0 \frac{\partial}{\partial t_i}\o(x,st)ds$ is
in $\O^\ell(\log D)$, for taking as equation for $D$ the equation
$h$ of $E$ (regarded as a function of $(x,t)$ independent of $t$) we
have 
$$h(x)\int^1_0\frac{\partial}{\partial t_i}\o(x,st)ds
=\int^1_0\frac{\partial}{\partial t_i}h(x)\o(x,st)ds$$
and is therefore regular, and similarly
$$dh\wedge\int_0^1\frac{\partial}{\partial t_i}\o(x,st)ds
=\int_0^1\frac{\partial}{\partial t_i}dh\wedge\o(x,st)ds,$$
also regular. 

Hence
$$\o\equiv \o(x,0)\quad (\mbox{mod}\ m_{U,0}\O^\ell(\log D))$$
and the required inclusion will follow by Nakayama's Lemma,
once we show that 
$$\o(x,0)\in p^\ast(\O^\b(\log E))\wedge \O^\b_U.$$
But this is now evident: if
$$\o(x,0)=\sum dt_{i_1}\wedge\cdots\wedge dt_{i_k}\wedge \o_{i_1,\cdots,i_k}$$
with each $\o_{i_1,\cdots,i_k}$ independent of $t$ and having no $dt_i$
component, then each $\o_{i_1,\cdots,i_k}$ is $D$-logarithmic (as it is
the contraction of $\o(x,0)$ by a $D$-logarithmic $k$-vector) and, since
$\o_{i_1,\cdots,i_k}=p^\ast q^\ast \o_{i_1,\cdots,i_k}$, we have 
$\o_{i_1,\cdots,i_k}\in p^\ast(\O^\b(\log E))$.
\eop
Part (2) of the proposition now follows: we have
$$
\ch^\b_{D/S}
=
\frac{\ch^\b_D}{\sum ds_i\wedge \ch^{\b-1}_D}
=
\frac{\O^\b_{V\times S}}{\langle i^\ast(h_E\O^\b(\log E))\rangle
+\sum ds_i\wedge\O^{\b-1}_{V\times S}}
=
\frac{\O^\b_{V\times S}}{h_D\O^\b(\log D)+\sum ds_i\wedge\O^{\b-1}_{V\times S}}.
$$
\eop
\begin{proposition}
If $E\subset W$ is a free divisor, $i_0:V\to W$ is a map and
$D_0=i_0^{-1}(E)$, then we have the following inclusions:
$$
h_{D_0}\O^k_V+dh_{D_0}\wedge\O^{k-1}_V
\subseteq
\langle i_0^\ast(h_E\O^\b(\log E)\rangle_k
\subseteq
h_{D_0}\O^k(\log D_0)
$$
where the lower index $k$ in the middle complex means the degree $k$
summand.
If $D_0$ is almost free, then 
provided $\dim\ D_0\geq 2$ and $E$ is not smooth, the first inclusion is strict.
If $i_0$ is not algebraically transverse to $E$ then
the second inclusion is also strict for $k= \dim\ D_0$ and $k=\dim\ D_0+1$.
\end{proposition}
\Proof 
The first inclusion holds because
$dh_E/h_E\in\O^1(\log E)$ and 
we can take $h_{D_0}=i_0^\ast(h_E)$; thus
$h_{D_0}\in i_0^\ast(h_E\O^0(\log E))$ and 
$dh_{D_0}\in i_0^\ast(h_E\O^1(\log E))$.
The second inclusion holds because, 
as is easily checked directly from the definition,
$i_0^\ast(\O^\b(\log E))\subset \O^\b(\log D_0)$.

Strictness of the first inclusion when $D_0$ is almost free holds by
\ref{con}; for this implies that the quotient of the right hand side by the
left is supported at all points of $D_{0,\mbox{sing}}\setminus\{0\}$,
and thus also at $0$.

Strictness of the second inclusion for $k=\dim\ D_0$ follows from
the fact that the torsion submodule of $\ch^p_{D_0}$ has length
equal to that of $N{\cal K}_{E,e}i$ (by \ref{torsion}, below), 
and thus $\depth\ \ch^p_{D_0}=0$. The depth of $\O^p_V/h_{D_0}\O^p(\log D_0)$,
on the other hand, is at least 1; for
$\O^p(\log D_0)$, and thus $h_{D_0}\O^p(\log D_0)$, are 
isomorphic to $\Der(\log D_0)$ and thus
have depth at least 2, from which it follows that 
$\O^p_V/h_{D_0}\O^p(\log D_0)$ has depth at least 1.

For $k=\dim\ V$, it is evident that
$i_0^\ast(h_E\O^k(\log E))\subseteq m_{V,0}\O^k_V$
unless $i_0$ is algebraically transverse to $E$.
\eop
\begin{proposition}\label{think} Let $D_0$ be an almost free divisor.
Then for each vector field $\chi\in\Der(\log D_0)$, contraction by
$\chi$ gives rise to a well-defined $\OO_{D_0}$-linear morphism
$\ch^k_{D_0}\to\ch^{k-1}_{D_0}$.
\end{proposition}
\Proof The content of the statement is that if we represent $\o\in\ch^k_{D_0}$
by $\o_1\in\O^k_V$ then the class of 
$\iota_{\chi}(\o_1)$ in $\ch^{k-1}_{D_0}$
is independent of the choice of $\o_1$.
Since contraction is evidently linear, this amounts to showing that
if $\o_1\in\O^k_V$ is $0$ in $\ch^k_{D_0}$ in  then $\iota_{\chi}(\o_1)$ is
$0$ in $\ch^{k-1}_{D_0}$. 

Let $i:V\times S\to W$ be a deformation of the inducing map $i_0$ which 
is algebraically
transverse to $E$, and let $D=i^{-1}(E)$. 
Denote by $j:V\to V\times S$ the
inclusion $j(v)=(v,0)$. We have $D_0=j^{-1}(D)$, and since $D$ is free we
can apply Proposition \ref{gasp} to conclude that 
$\Der(\log D_0)=tj^{-1}(j^\ast
(\Der(\log D)))$. As $j$ is an inclusion this means simply that
for each $\chi\in\Der(\log D_0)$ there exists $\chi_2\in \Der(\log D)$
extending $\chi$. In particular, if $\chi_2=\sum a_i\partial /\partial v_i
+\sum b_j\partial /\partial s_j$ then all $b_j$ vanish when $s=0$.

Let $\o_2\in\O^k_{V\times S}$ extend $\o_1$. As
$\o_1$ is equal to $0$ in $\ch^k_{D_0}$, we have
$$
\o_2\in h\O^k(\log D)+\sum s_i
\O^k_{V\times S}+\sum ds_i\wedge \O^{k-1}_{V\times S}
.$$
Now let $\chi_2\in\Der(\log D)$ extend $\chi$. Then 
$\iota_{\chi}(\o_1)=j^\ast(\iota_{\chi_2}(\o_2))$. Since
$$\iota_{\chi_2}(h\O^k(\log D))\subset h\O^{k-1}(\log D),$$
$$\iota_{\chi_2}(\sum s_i\O^k_{V\times S})\subset \sum s_i \O^{k-1}_{V\times S}
$$
and
$$\iota_{\chi_2}(\sum ds_i\wedge \O^k_{V\times S})\subset 
\sum ds_i\wedge\O^{k-2}_{V\times S}+\sum s_i \O^{k-1}_{V\times S},$$
it follows that 
$\iota_{\chi}(\o_1)=j^\ast(\iota_{\chi_2}(\o_2)
\in j^{\ast}(h\O^{k-1}(\log D)+\sum s_i\O^{k-1}_{V\times S}+\sum ds_i\wedge
\O^{k-2}_{V\times S})$
and the proposition  is proved.
\eop
We end this section with a definition of logarithmic critical space:
\begin{definition}  Let $\pi:D\to S$ be a deformation of an almost
free divisor and choose a free extension 
$\rho:{\cal D}\to S\times T$ of $\pi$.
Then 
\be 
\item
the {\em logarithmic critical space} of $\rho$, $C^{\log}_\rho$, 
is $$C^{\log}_\rho=\mbox{supp}T^{1,\log}_{{\cal D}/S\x T}$$
with analytic structure defined by the ideal
${\cal C}^{\log}_{\rho}:=
{\cal F}_0^{\OO_{V\times S\times T}}(T^{1,\log}_{{\cal D}/S\x T})$
\item 
the logarithmic critical space of $\pi$ is 
$C^{\log}_{\pi}= C^{\log}_\rho\cap \rho^{-1}(S\times\{0\})$, 
with its natural analytic 
structure as an intersection.(Once again,  
this does
not depend on the choice of free extension $\rho$ of $\pi$.) 
\ee
\end{definition}
Here ${\cal F}_0$ is the zero'th Fitting ideal. We remark that in case
the codimension of $C^{\log}_p$ in $V\times S\times T$ is equal to 
$\dim S\times T
+1$, then this coincides with the annihilator of the 
${\OO_{V\times S\times T}}$-module $\t(\pi)/t\pi(\Der(\log D)$, by the
theorem of \cite{waam}. This is the case, for example, if 
$\pi:D\to S$ is a free deformation of
an almost free divisor $D_0$, and frees $D_0$.
The  ideal ${\cal C}^{\log}_\pi$
can be calculated very easily: if $\Theta$ is a
matrix whose columns are the components of a free basis of
$\Der(\log {\cal D})$
and $\Theta'$ its submatrix consisting of the rows corresponding
to the deformation parameters, then
${\cal C}^{\log}_{\rho}$
is generated by the restriction to $V\x S$ of the maximal minors of $\Theta'$.

One of the major differences with the classical case is that 
${\cal C}^{\log}_\pi$
is not necessarily radical even when $\pi:D\to S$ is the deformation
induced by a ${\cal K}_{E,e}$-versal deformation of the inducing map
$i$
(cf \cite{damonlegacy}). \\
\bc
\section{Algebraic properties}
\label{algebraic properties}
\ec
As remarked above, our complex $\check{\O}^\b_{D/S}$ has better depth 
properties than the standard complex $\O^\b_{D/S}$. The improvement 
is crucial. 
\begin{lemma}\label{depth}
Let $D\subseteq\CC^N$ be a free divisor, and let $g_1,\cdots,g_s$ and
$f_1,\cdots,f_m$ be holomorphic functions on $D$, such that
 $g_1,\cdots,g_s$ is 
a regular sequence on $\OO_D$, and such that
for each $i$, $\dim\ C^{\log}_g\cap V(g_1,\cdots,g_i)
< \dim\ D\cap V(g_1,\cdots,g_i)$ (where $C^{\log}_g$ is the intersection
with $D$ of the support of $\theta(g)/tg(\Der(\log G))$). 
Denote by $f$ the map
$(f_1,\cdots,f_m):D\to\CC^m$. 
Let $\OO=\OO_D/(g_1,\cdots,g_s)$ and
$$\ch^k:=\ch^k_D/(g_1,\cdots,g_s)\ch^k_D.$$
$$\ch^k(i)=\ch^k/\sum_{j=1}^idf_j\wedge\ch^{k-1}.$$
Then\\ 

(1) for $0\leq i\leq m$, and $0\leq k < \dim\ V(g) - \dim\ C^{\log}_f
\cap V(g)$, $\mbox{depth}\ \ch^k_x(i)\geq \dim V(g) -k$, and\\

(2) the map $$df=df_1\wedge\cdots\wedge df_m:\frac{\ch^k_{D,x}}
{(g_1,\cdots,g_s)\ch^k_{D,x}
+\sum_{j=1}^s df_j\wedge\ \ch^{k-1}_{D,x}}\to \ch^{k+m}_x$$
is injective, for $0\leq k <\mbox{depth}\ \OO_D/(g_1,\cdots,g_s) -
\dim_x C^{\log}_f\cap V(g)$.
\end{lemma}
\Proof This lemma is just a translation of Lemma 1.6 of \cite{greuel}
(which can also be found in section 8C of \cite{looijenga})
and the proof is essentially the same. First, note that because
$\ch^k_D$ is a maximal Cohen-macaulay $\OO_D$-module, any $\OO_D$-regular
sequence is also $\ch^k_D$-regular, and thus $\ch^k$ is a maximal
Cohen-Macaulay $\OO$-module. 
We now follow Greuel's argument. Write $f(x)=z$, and $f^{-1}(z)=D_z$.\\

(a) If $x\notin C^{\log}_f\cap V(g)$ then $D,x$ is isomorphic over
$\CC^s,z$ to $D_z\times \CC^s$. The map
$$df_1\wedge \cdots \wedge df_s\wedge:\ch^k_{D/\CC^s}\to \ch^{k+s}_D$$
then has a left inverse, and it follows that it remains exact after
tensoring with $\OO_x$.\\

(b) Suppose that $x\in C^{\log}_f\cap V(g)$. 

We prove (1) by induction on $k$. For $k=0$ it is trivial. Suppose
it true for $k-1$, and consider the sequence 
$$0\to \ch^{k-1}_x(i)\stackrel{df_i\wedge}{\longrightarrow}
\ch^k_x(i-1)\to\ch^k_x(i)\to 0,$$
which is evidently exact except perhaps at $\ch^{k-1}_x(i)$. To see that
$df_i\wedge$ is injective, note that it is certainly injective off
$C^{\log}_f\cap V(g)$, and thus $\ker df_i\wedge\subseteq{\cal H}^0_{C^{\log}_f
\cap V(g)}(\ch^{k-1}_x(i))$. However, $\depth\ch^{k-1}_x(i) > \depth
\OO_x-k > \dim C^{\log}_f \cap V(g)$, and so this local cohomology group
vanishes and the sequence is exact.
Now we argue by induction on $i$: for $i=0$, $\depth\ \ch^k_x(i)
=\depth\ \ch^k_x=\depth\ \OO_x$, as shown at the start of this proof. 
Suppose that $\depth\ \ch^k_x(i-1)\geq \depth\ \OO_x-k$. Then by exactness
of the above sequence, it follows that $\depth\ \ch^k_x(i)\geq \depth\ \OO_x-k$,
as required. This completes the proof of (1). 

Part (2) follows, for now
by (a), $\ker df\wedge:\ch^k_x/\sum_{j=1}^m df_j\wedge \ch^{k-1}_x\to
\ch^{k+m}_x$ is contained in ${\cal H}^0_{C^{\log}_f \cap V(g)}(
\ch^k_x/\sum_{j=1}^m df_j\wedge \ch^{k-1}_x)$, which, by the depth estimate
just obtained, is equal to $0$.\eop

We now apply the lemma in the following situation: $\pi:D\to S$
is an admissible deformation of an almost free divisor, 
$\rho:{\cal D}\to S\x T$ 
is a free extension of 
$\pi$, and denoting coordinates on $S$ and $T$ by $s_1,\cdots,s_d$ and
$t_1,\cdots,t_e$,  we take 
$(g_1,\cdots,g_k)=(t_1,\cdots,t_e)$ and $(f_1,\cdots,f_m)
=(t_1,\cdots,t_e,s_1,\cdots,s_d)$.  From the lemma we conclude immediately
that
\begin{proposition}\label{block} If $\pi:D\to S$ is an admissible deformation
of $D_0$, and $\rho:{\cal D}\to S\x T$ is a free extension of 
$\pi$, then 

(1) provided $0\leq k < \codim\ C^{\log}_{\pi}$, 
$\depth\ \ch^k_{D/S}\geq \dim D_0 - k$ 
and 
$d\pi:\ch^k_{D/S}\to\ch^{k+d}_D$ is an injection.\\

(2) If $D_0$ is an almost free divisor and $\pi:D\to S$ is a deformation 
which frees $D_0$,
then
$\codim\ C^{\log}_{\pi} =\dim\ S+1$, so that the assertions of (1) hold for
$k\leq \dim D_0$.\\

(3) The map 
$\lambda d\pi:\ch^k_{D/S}\to
\ch^{k+e}_{\cal D}/(t_1,\cdots,t_e)\ch^{k+e}_{\cal D}$
defined by wedging with $ds_1\wedge\cdots\wedge ds_d\wedge 
dt_1\wedge\cdots\wedge dt_e$, is 1-1
for $0\leq k < \codim\ C^{\log}_{\pi}$.\eop
\end{proposition}

\begin{proposition}\label{torsion} Let $\pi:D\to S$ be an admissible
deformation of the $p$-dimensional divisor $D_0$.
Then
\be
\item
$\ch^k_{D/S}$ is a torsion-free $\OO_D$-module for 
$0\leq k <\codim\ C^{\log}_\pi$ (in the sense that every non-zero-divisor in
$\OO_D$ is regular on $\ch^k_{D/S}$);
\item
$\ch^k_{D/S}$ is a torsion module for $k = \dim\ D_0+1$; (and recall that
$\ch^k_{D/S}=0$ for $k>\dim\ D_0+1$).
\item
if $\rho:{\cal D}\to S\x T$ is a free extension of $\pi$ then
the cokernel of $\lambda d\pi:\ch^p_{D/S}\to
\ch^{p+d+e}_{\cal D}\otimes \OO_D$ is isomorphic to $T^{1,\log}_{D/S}$.
This cokernel is independent of choice of $\rho$;
\item if $D_0$ is almost free then the torsion submodule of
$\ch^{p}_{D_0}$ has dimension equal to 
$\dim T^{1,\log}_{D_0}$.
\ee
\end{proposition}
\Proof 1. Let $p:{\cal D}\to S\x T$ be a free extension of $\pi$, 
with 
parameter space $T$ of dimension $e$. 
Write $d=\dim\ S$. For $0\leq k <\codim\ C^{\log}_{\pi}$, 
$\lambda d\pi:\ch^k_{D/S}\to\ch^{k+d}_{D}$ is an injection, 
so it remains
only to show that $\ch^{k+d}_D$ is torsion-free. This follows by
the last part of Proposition~\ref{block}: $dh\wedge$ 
embeds $\ch^{k+d}_D$ in the free $\OO_D$-module $\O^{k+d+1}_{V\x S}/
h\O^{k+d+1}_{V\x S}$.\\

2. This is clear, since if $k > \dim D_0$ then 
$\ch^k_{D/S}$ is supported on $C^{\log}_{\pi}$, whose 
dimension is less than that of $D$.\\

3. Contraction of the generator $dy\wedge ds\wedge dt$ of 
$\O^{p+d+e+1}_{V\times S\x T}$ by vector fields gives an isomorphism 
$\iota:\t_{V\x S\x T}\to \O^{p+d+e}_{V\x S\x T}$ 
which restricts to an isomorphism
$\Der(\log {\cal D})\simeq h\O^{p+d+e}(\log {\cal D})$. Thus 
$$\frac{\t_{V\x S\x T}}{\Der(\log \D)}\simeq \ch^{p+d+e}_{\D}.$$
The preimage under $\iota$ of the submodule $dt\wedge 
ds\wedge\O^{p}_{V\times S\x T}$
of $\O^{p+d+e}_{V\x S\x T}$ 
is the module of vector fields with no component in the $S\x T$ direction,
i.e. $\ker t\rho:\t_{V\x S\x T}\to\t(\rho)$; putting everything
together, it follows that 
$$\frac{\ch^{p+d+e}_{\D}}{dt\wedge ds\wedge\ch^{p}_{D/S}}\otimes\OO_D
\simeq \frac{\t_{V\x S\x T}}{\Der(\log \D)+\ker(t\rho)}\otimes \OO_D
\simeq \frac{\t(\rho)}
{t\rho(\Der(\log \D)}\otimes \OO_D\simeq T^{1,\log}_{D/S}.$$
Independence of the choice of $\rho$ follows as usual from 
Lemma~\ref{welldef}.

4. Let $\pi:D\to S$ be a free deformation of $D_0$, and let $\dim\ S = d$. 
Then
the torsion submodule of $\ch^{p}_{D_0}$ is equal to the kernel of
the map 
$$ds_1\wedge\cdots\wedge ds_d:\ch^{p}_{D_0}
=\frac{\ch^{p}_{D/S}}{(s_1,\cdots,s_d)\ch^{p}_{D/S}}\to
\frac{\ch^{p+d}_{D}}{(s_1,\cdots,s_d)\ch^{p+d}_{D}},$$
for composition with the injective morphism $dh\wedge$ maps $\ch^{p}_{D_0}$
into the free $\OO_{D_0}$-module\\ 
$\O^{p+d+1}_{V\x S}/(h, s_1,\cdots,s_d)
\O^{p+d+1} _{V\x S}$ (so that $T\ch^{p}_{D_0}$ is 
contained in the kernel), while on the other hand the kernel is supported 
only at $0$, since $D_0$ is free elsewehere, and hence is a torsion module.
Now this kernel can be identified with 
$\Tor_1^{\OO_S}(T^{1,\log}_{D/S},\OO_S/m_{S,0})$; 
for tensoring the short exact sequence
$$0\to \ch^p_{D/S}\stackrel{ds\wedge}{\too}\ch^{p+d}_D
\to T^{1,\log}_D/S\to 0$$
with $\CC=\OO_S/m_{S,0}$ we get the long exact sequence
$$0=\Tor_1^{\OO_S}(\ch^{p+d}_D,\CC)\to \Tor_1^{\OO_S}(T^{1,\log}_{D/S},\CC)
\to \ch^p_{D_0}\to \ch^{p+d}_D\otimes\CC\to T^{1,\log}_{D_0}\to 0.$$

Now assume furthermore that $\pi:D\to S$ is versal (i.e. the inducing map
$i:V\times S\to W$ is a ${\cal K}_{E,e}$-versal deformation of the
map $i_0$ which induces $D_0$). Then there is
an exact sequence
\begin{equation}
0\to {\cal L}\stackrel{\alpha}{\too}
\t_S\stackrel{\rho}{\too}\pi_\ast(T^{1,\log}_{D/S})\to 0 
\end{equation}
in which $\rho$ is the Kodaira-Spencer map of the deformation. As
 $T^{1,\log}_{D/S}$ is a Cohen-Macaulay $\OO_{V\x S}$-module of dimension
$d-1$, by~\ref{cm}, and 
is finite
over $\OO_S$, it follows that its $\OO_S$-depth is equal to $d-1$ also, and it
follows that $\cal L$ is free, of rank $d$. 
If $\pi$ is miniversal, then all entries
in the matrix $\alpha$ lie in $m_{S,0}$ and $\alpha\otimes \CC=0$. Hence
$\Tor^{\OO_S}_1(T^{1,\log}_{D/S},\CC)\simeq {\cal L}/m_{S,0}{\cal L}$ has
dimension $d=\dim\ S = \dim\ T^{1,\log}_{D_0}$.
\eop

Part (4) of this result, together with Damon's theorem identifying
the ${\cal A}_e$-normal space of a map-germ $f_0:\CC^n,0\to\CC^p,0$
with $T^{1,\log}_{D(f_0)}$ (re-proved in Section \ref{last} below), 
now proves Theorem \ref{tor}.
\begin{example}
{\em
Suppose that  $E\subset W$ is a weighted homogeneous free divisor 
and $i:V\to W$ is
weighted homogeneous of ${\cal K}_{E,e}$-codimension 1. Let $D_0=i^{-1}(E)$,
let $x_1,\cdots,x_p$
be coordinates on $V$, and denote by $\chi_e$ the Euler vector field on $V$.
Then the torsion in $\ch^p_{D_0}$ is generated over $\CC$ by the class in 
$\ch^p_{D_0}$ of  
$\iota_{\chi_e}(dx_1\wedge\cdots\wedge dx_p).$
For $\ch^p_{D_0}$ is supported precisely at those points where $i$ is not 
algebraically transverse to $E$, and thus $dx_1\wedge\cdots\wedge dx_p$
is a torsion form in $\ch^p_{D_0}$. As the Euler field is logarithmic,
$\iota_{\chi_e}(dx_1\wedge\cdots\wedge dx_p)$ is also a torsion form
(the contraction is well-defined in $\ch^{p-1}_{D_0}$ by \ref{think}).
It is not zero, since its exterior derivative is a non-zero multiple of 
$dx_1\wedge\cdots\wedge dx_p$.
}
\end{example}
\section{Using $\ch^\b_{D/S}$ to calculate cohomology}\label{using}
\begin{definition} The analytic space $D\subseteq W$ is {\em locally
quasihomogeneous} if for each $x\in D$ there exist local analytic
coordinates for $W$, centred at $x$, with respect to which the germ $D,x$
is weighted homogeneous (with strictly positive weights).
\end{definition}
This property played an important r\^{o}le in \cite{CMN} (under the
name {\em strong quasihomogeneity}):
\begin{theorem}\label{CMN} {\em \cite{CMN}} Let $D\subset \CC^n$ be a locally
quasihomogeneous free divisor, and let $U=\CC^n\setminus D$. Then integration
along cycles induces an isomorphism
$$h^q(\Gamma(\CC^n,\O^\b(\log D)))\simeq H^q(U;\CC)$$
for $0\leq q\leq n$.
\eop
\end{theorem}
Here and in what follows, we use the small $h^q$ for the $q$'th
homology of a complex.

One can speak of a locally quasihomogeneous germ - one which
has a locally quasihomogeneous representative. This allows the 
apparent oxymoron
of a quasihomogeneous germ which is not locally 
quasihomogeneous; however this incongruity
seems less disturbing
than the possibility of a strongly
quasihomogeneous space which is not (globally) quasihomogeneous, so the term
``locally quasihomogeneous'' is to be preferred to 
``strongly quasihomogeneous''. 

By (~\ref{absbits}).6, if $D\subseteq V$ is a locally weighted homogeneous
divisor then
the complex 
$$0\to\CC_D\to \OO_D\to\ch^1_D\to\ch^2_D\to\cdots$$
is a resolution of $\CC_D$. 
It follows that if $V$ is a Stein space then 
$$h^q(\G(D, \ch^\b_D))\simeq H^q(D;\CC).$$
Hence we can use the complex $\ch^k_{D/S}$ to calculate the cohomology of
the free fibres of an admissible deformation, provided the free fibres are
locally quasihomogeneous.
\begin{example} {\em 1. If $(\dim\ U,\dim\ V)$ are in Mather's range of
nice dimensions (cf \cite{mather6}) 
then every stable map-germ $U\to V$ is quasihomogeneous
with respect to some coordinate system. (In fact this characterises the nice 
dimensions). It follows that in the nice
dimensions stable discriminants are locally quasihomogeneous. Thus, for
our motivating example, the insistence on local quasihomogeneity is not
unduly restrictive.}
\end{example}

\begin{definition} Let $E\subseteq W$ be a divisor and let $\cal S$
be the canonical Whitney stratification of $E$. 
\be
\item
A stratum in $\cal S$
is {\em locally quasihomogeneous} if at each point $x\in S$, there are local
analytic coordinates centred at $x$ with respect to which $E$ is weighted
homogeneous (with respect to strictly positive weights). The
{\em weighted homogeneous codimension} $\mbox{wh}(E)$ 
of $E$ is the infimum of the codimensions (in $E$) of strata of $\cal S$
which are not
locally quasihomogeneous, and 
$\infty$ 
if $E$ is locally quasihomogeneous. 
\item
A stratum $S\in {\cal S}$ is {\em holonomic} if at each point $x\in S$,
$T_xS=T^{\log}_xE$. The {\em holonomic codimension} $\mbox{hn}(E)$of $E$ is 
the infimum of the codimensions (in $E$) of strata of $\cal S$ which are not
holonomic, and 
$\infty$       
if every stratum is holonomic.\ee
\end{definition}
\begin{lemma}\label{res} Let $D_0$ be an almost free divisor 
based on $E$, suppose $\dim\ D_0\leq \mbox{wh}(E)$, and 
let $\pi:D\to S$ be an admissible deformation of $D_0$.
Then for
$x\notin C^{\log}_\pi$, the complex $\ch^\b_{D/S}$ is 
a resolution of $\pi^{-1}(\OO_S)$. 
\end{lemma}
\Proof Let $i:V\x S\to W$ induce $D$ from $E$ (i.e. $D=i^{-1}(E)$.) 
and for $s\in S$ let $i_s:V\to W$ be defined by $i_s(x)=i(x,s)$.
By a variant of the standard transversality lemma, 
if $(x,s)\notin C^{\log}_\pi$ the
$i_s$ is algebraically transverse to $E$ at $x$, and thus transverse to
each of the strata in the canonical Whitney stratification of $E$. As
$\dim V < \mbox{wh}(E)$, $i_s$ meets only strata along which $E$ is locally
weighted homogeneous, and thus $D_s=i_s^{-1}(E)$ is locally weighted
homogeneous. It follows that $\ch^\b_{D_s,x}$ is a resolution of $\CC_{D_s,x}$.

Again because $(x,s)\notin C^{\log}_\pi$,  there is a commutative diagram
$$\begin{array}{ccccc}
(D,(x,s))&&\stackrel{\simeq}{\too}&&(D_s,x)\x (S,s)\\
&\pi\searrow&&\swarrow\mbox{projection}&\\
&&(S,s)&&
\end{array}$$
and the isomorphism $(D,(x,s))\simeq (D_s,x)\x (S,s)$ induces an isomorphism
of complexes 
$$\ch^\b_{D/S}\simeq \ch^\b_{D_s}\otimes_{\CC} \pi^{-1}(\OO_S)$$
and thus $\ch^\b_{D/S}$ is a resolution of $\CC\otimes\pi^{-1}(\OO_S)=
\pi^{-1}(\OO_S)$.
\eop

Our notion of {good representative} of a germ of
admissible deformation $\pi:(D,(x_0,0))\to S,0$ of an almost free
divisor $D_0,x_0$ is an obvious modification
of the notion standard in the theory of isolated singularities, as described
in some detail in \cite{looijenga}, pp 25-26. We require  
an open set $X\subset \CC^{p+1}\x \CC^d$, a map $i:X\to 
W$ such that, writing $\tilde{D} = i^{-1}(E)$, we have 
$\tilde{D}\cap \CC^{p+1}\x \{0\}=D$ (as germs at $x_0$), 
a real-analytic function 
$r:X\to [0,\infty)$  and a real $\varepsilon > 0$ such that 
$r^{-1}(0)\cap D_0 = \{x_0\}$ and such that 
$\tilde{D}_0$ is stratified transverse to $r^{-1}(\varepsilon')$ (in
$\CC^{p+1}$)  for each $\varepsilon'$ with $0 < \varepsilon'\leq \varepsilon$.
Then there exists a contractible $S\subset \CC^d$ containing $0$ such
that the restriction of $\pi$ to each stratum of the canonical 
Whitney stratification of $\tilde{D}$  is a submersion
at each point of $\pi^{-1}(S)\cap r^{-1}(\varepsilon)$. Finally, we take
$A=r^{-1}([0,\varepsilon)) \cap \pi^{-1}(S)$, and $D=\tilde{D}\cap A$. \\

From now on we suppose that we are given a good Stein representative of
an admissible deformation $\pi:D,(x_0,0)\to S$ of the $p$-dimensional 
almost free divisor
$(D_0,x_0)$ 
based on the free divisor $E$. We suppose that $p < \mbox{min}\{\mbox{wh}(E),
\mbox{hn}(E)\}$, and denote the dimension of $D$ by $m$.
Consider the short exact
sequence of complexes
\begin{equation}\label{ses0}
0\to\H^0(\ch^\b_{D/S})\to\ch^\b_{D/S}\to\ch^\b_{D/S, (>0)}
\to 0
\end{equation}
where $\H^0(\ch^\b_{D/S})$ is the 0'th cohomology sheaf (i.e. 
$\ker(d:\ch^0_{D/S}\to\ch^1_{D/S})$) and $\ch^\b_{D/S, (>0)}$ is the complex
$0\to \ch^0_{D/S}/\ker (d:\ch^0_{D/S}\to\ch^1_{D/S}) \to \ch^1_{D/S}\to\cdots.
$
Applying the functor $\pi_\ast$ to this sequence and deriving,
we obtain
a long exact sequence of cohomology on $S$:
\begin{equation}\label{seq0}
\cdots\to R^q\pi_\ast(\H^0(\ch^\b_{D/S}))\to \RR^q\pi_\ast(\ch^k_{D/S})
\to \RR^q\pi_\ast(\ch^\b_{D/S, (>0)})\to\R^{q+1}\pi_\ast(\H^0(\ch^\b_{D/S}))
\to\cdots
\end{equation} 
for $q\geq 1$.

\begin{lemma} In these circumstances,\\

\ni (1) provided $\dim V > 1$, $\H^0(\ch^\b_{D/S})=\pi^{-1}(\OO_S)$, and so 
$R^q\pi_\ast(\H^0(\ch^\b_{D/S}))$ is canonically isomorphic to 
$R^q\pi_\ast(\CC_D)\otimes_{\CC}\OO_S$.

(2)
$$R^q\pi_\ast\ch^\b_{D/S, (>0)}\simeq \pi_\ast(
\H^q(\ch^\b_{D/S, (>0)})=\pi_\ast(\H^q(\ch^\b_{D/S})$$ for
$q\geq 1$.\\

(3)
$$R^q\pi_\ast(\ch^\b_{D/S})=\H^q(\pi_\ast\ch^\b_{D/S}).$$
\end{lemma}
\Proof (1) The statement is clearly true outside of $C^{\log}_\pi$,
and moreover it is obvious that
$\H^0(\ch^\b_{D/S})\supseteq \pi^{-1}(\OO_S)$ everywhere. The opposite
inclusion follows from the connectedness of the fibres of $\pi$.\\

\ni (2) 
As $D_0$ is an almost free divisor, it follows that $\pi_{|C^{\log}_\pi}$
is finite. For $q>0$, ${\cal H}^q(\ch^\b_{D/S})$ is supported only on
$C^{\log}_\pi$, and thus $\RR^i\pi_\ast {\cal H}^q(\ch^\b_{D/S})=0$ for
$i>0$. Thus the spectral sequence $\RR^p\pi_\ast
{\cal H}^q(\ch^\b_{D/S, (>0)})
\implies R^{p+q}\pi_\ast\ch^\b_{D/S, (>0)}$ collapses at $E_1$ and
the conclusion follows.\\

\ni 3. This is an immediate consequence of the $\OO_D$-coherence of the sheaves
$\ch^k_{D/S}$, together with the fact that we have chosen a Stein 
representative.
\eop

\begin{corollary}\label{les} Provided $\dim\ V > 1$, the long exact sequence 
(\ref{seq0}) reduces to
\begin{equation}\label{seq1}
\cdots \to R^q\pi_\ast(\CC_D)\otimes\OO_S\to \H^q(\pi_\ast(\ch^k_{D/S}))
\to \pi_\ast(\H^q(\ch^\b_{D/S}))\to\R^{q+1}\pi_\ast\CC_D\otimes\OO_S
\to\cdots
\end{equation}
for $q > 1$.
\eop
\end{corollary}

\begin{theorem}\label{coherence} Let $\pi:D\to S$ be a good representative of
an admissible deformation of the $p$-dimensional almost free divisor 
$D_0$, and suppose that
$p<\mbox{hn}(E)$. Then
${\cal H}^q(\pi_\ast(\ch^k_{D/S}))$ is a coherent
sheaf of $\OO_S$-modules.
\end{theorem}
\Proof We use a theorem of Duco van Straten (\cite{straten}).
Let $\pi:D\to S$ be a good
representative of an admissible deformation of an almost free divisor,
and let ${\cal L}$ be a sheaf of complex vector spaces on $D$. $\cal L$ is
{\em transversally constant} at the boundary if there exist an open
neighbourhood $U$ of $\partial A$ and a $C^\infty$ vector field 
$\t$ on $U$ with the following properties:
\be
\item
$\t$ is transverse to $\p A$;
\item
$\t$ is tangent to $D$ and to the fibres of $\pi$;
\item
the restriction of $\cal L$ to each local integral curve of $\t$ is a 
constant sheaf.
\ee
Then van Straten proves
\begin{theorem}\label{duco}
Let $\pi:D\to S$ be a good representative of a germ $\pi:(D,0)\to (S,0)$, and
let $K^\b$ be a finite complex of sheaves on $D$. Assume
\be
\item
the sheaves $K^i$ are $\OO_D$-coherent;
\item
the differentials $K^i\to K^{i+1}$ are $\pi^{-1}(\OO_S)$-linear;
\item
the cohomology sheaves ${\cal H}^i(K^\b)$ are transversally constant at the 
boundary. 
\ee
Then $\RR^i\pi_\ast K^\b$ is an $\OO_S$-coherent module.
\eop
\end{theorem}
The statement in \cite{straten} assumes a $1$-dimensional base $S$, but this 
hypothesis is not used in the proof.\\

\begin{lemma} Let $\pi:D\to S$ be a good representative of an admissible
deformation of the $p$-dimensional  
almost free divisor $D_0$ based on the free divisor $E$, 
with $p < \mbox{hn}(E)$. Then the cohomology sheaves
${\cal H}^i(\ch^\b_{D/S})$ are transversally constant at the boundary.
\end{lemma}
\Proof 
As $C^{\log}_\pi\cap D_0=\{(0,0)\}$, at each point of the boundary
$\p A\cap D_0$ all the strata of the canonical 
Whitney stratification ${\cal X}_0$
of $D_0$ are holonomic, i.e. for each stratum $X\in {\cal X}_0$ and $x\in\p D$,
$T_xX=T^{\log}_xD_0$. Thus,
since $\p A_0$ is transverse to ${\cal X}_0$, at each point $x\in \p D_0$ 
there exists a germ of vector field $\chi_x^0\in\Der(\log D_0)_x$ such that
$\chi_x^0(x)\notin T_x\p A$. As $\pi$ is locally trivial at $(x,0)$,
the logarithmic germ $\chi_x^0$ extends to a germ 
$\chi_x\in \Der(\log D)_{(x,0)}\cap \ker (d\pi)$; there is a neighbourhood
$U_x$ of $(x,0)$ in $V\x S$ such that for all $(x',s)\in \p D\cap U_x$, 
$\chi_x(x',s) \notin T_{(x,s)}\p A$. Let $U=\bigcup_{x\in D_0\cap\p A}
U_x$.  Now by means of a partition of unity subordinate to the open 
cover $\{U_x\}$ of $U$, we construct from these $\chi_x$ a $C^\infty$
vector field $\t$ on $U$ such that 
\be
\item
for all $(x,s)\in U$, $\t(x,s)\in T^{\log}_{(x,s)}D$
\item
for all $(x,s)\in U\cap\p D$, $\t(x,s)\notin T_{(x,s)}\p A$. 
\ee
By shrinking $U$ if necessary, we can suppose that $U\cap C^{\log}_\pi =
\emptyset$. Thus at each point $(x,s)\in U$, $\pi$ is locally
analytically trivial. As $\t$ is tangent to the fibres of $\pi$ and 
logarithmic, it follows that the sheaves ${\cal H}^i(\ch^\b_{D/S})$ are
constant along the integral curves of $\t$.
\eop
Coherence of the sheaves $\RR^q\pi_\ast\ch^\b_{D/S}$ now follows, by
van Straten's theorem.
\eop

As in the classical case, that is, the theory of isolated complete 
intersection singularities (ICIS), 
from the coherence theorem follow a number of
interesting results. In most cases proofs are identical with those of the
corresponding classical results in 
\cite{looijenga}, and are omitted where possible.

Following Brieskorn and Greuel, we define three associated modules.
Let $\pi:D\to S$ be a good representative of an admissible deformation,
and let $\rho:{\cal D}\to S\x T$ be a free extension. We write
$\dim\ S = d,\  \dim\ T = r$.
For brevity let us denote by $\check{\cal H}$ the module ${\cal H}^p(\pi_\ast
\ch^\b_{D/S})$.
$$\check{\cal H}'= \frac{\pi_\ast\ch^p_{D/S}}{d(\pi_\ast \ch^{p-1}_{D/S})};$$
$$\check{\cal H}''= \frac{\pi_\ast\ch^m_D}{ds\wedge 
d(\pi_\ast \ch^{p-1}_{D/S})}.$$
$$\check{\cal H}''' = \frac{\pi_\ast (\ch^{p+d+r}_{\cal D})}
{\lambda ds\wedge d(\pi_\ast \ch^{p-1}_{D/S})}.$$

Note that by Lemma~\ref{welldef}, $\chch '''$ does not depend on the choice
of free extension $\rho$ of $\pi$. We also define
$$\check{H}= {\cal H}^p(\ch^\b_{D/S,(x_0,0)});$$
$$\check{H}'= \frac{\ch^p_{D/S,(x_0,0)}}{d \ch^{p-1}_{D/S,(x_0,0)}};$$
$$\check{H}''= \frac{\ch^m_{D,(x_0,0)}}{d\pi\wedge
d \ch^{p-1}_{D/S,(x_0,0)}}.$$
$$\check{H}''' = \frac{\ch^{p+d+r}_{{\cal D},(x_0,0)}}
{\lambda ds\wedge d \ch^{p-1}_{D/S,(x_0,0)}}.$$
By a standard argument (see e.g \cite{looijenga} 8.6), 
if $\pi:D\to S$ is a good representative of 
$\pi:(D,(x_0,0))\to S,0$ then the stalks at $0\in S$ of 
$\chch,\chch ',\chch ''$
and $\chch '''$ are canonically
isomorphic to $\chh, \chh ', \chh ''$ and $\chh '''$ respectively.

\begin{proposition} Let $D_0,x_0$ be a $p$-dimensional almost free divisor
based on the free divisor $E$ with $p < \mbox{hn}(E)$, and let  
$\pi:D\to S$ be a good Stein representative of an 
admissible
deformation of $D_0,x_0$. Then 
$\chch ', \chch ''$ and $\chch '''$
are $\OO_S$-coherent. 
\end{proposition}
\Proof  Apply the argument of~\ref{coherence} to the complexes
$\ch^\b_{D/S,(\leq p)}, d\pi(\ch^\b_{D/S,(\leq p)})$ and 
$\lambda d\pi(\ch^\b_{D/S,(\leq p)})$.\\
\eop

\begin{proposition}\label{exacta}
Let $D,x_0$ be a $p$-dimensional almost free divisor based on the free
divisor $E$, with $\mbox{wh(E)}>p+1$. Then
the complex 
$$0\to\CC_{D,x_0}\to\OO_{D,x_0}\to\ch^1_{D,x_0}\to \cdots \to \ch^p_{D,x_0}$$
is exact, $H^0_{\{x_0\}}(d\ch^{p-1}_{D,x_0})=0$, and there is a natural
chain of injections
$$H^1_{\{x_0\}}(d\ch^{p-2}_D)\hook H^2_{\{x_0\}}(d\ch^{p-3}_D)
\hook
\cdots \hook H^{p-1}_{\{x_0\}}(d\OO_D)\hook H^p_{\{x_0\}}(\CC_D).$$
\end{proposition}
\Proof Looijenga's proof of the corresponding result (\cite{looijenga}8.19) 
in the case
of a deformation $f:X\to S$ of an ICIS applies without change in this context. 
It relies only on estimates of
the depth of the modules $\OO^q_{X,x_0}$, 
and on the acyclicity of $\O^\b_{X,x_0}$
off $x_0$. 
Our depth estimates here (\ref{block}) are formally identical
to the estimates in the classical case (compare 8.15 and 8.16 of 
\cite{looijenga}), and the hypothesis that $p+1 <\mbox{wh}(E)$ guarantees
the acyclicity of $\ch^\b_{D}$ off $C^{\log}_\pi$. 
\eop
From this result follows a similar statement in the relative case:
\begin{proposition}\label{exactr}
Suppose that $\pi:D\to S$ is an admissible deformation of the $p$-dimensional
almost free divisor $D_0,x_0$ based on the free divisor $E$, 
and that $p+1<\min\{\mbox{wh}(E),\mbox{hn}(E)\}$. 
Then $\chh '$,  and $\chh '''$ are free $\OO_{S,0}$
modules of rank $\mu_\Delta(D_0,x_0)$, and if $\pi$ is a deformation
which frees $D_0$, then
the complex 
$$0\to\OO_{S,0}\stackrel{f^\ast}{\too}\OO_{D/S,(x_0,0)}\to\ch^1_{D/S, (x_0,0)}
\to\cdots\to\ch^p_{D/S,(x_0,0)}$$
is exact.
\end{proposition}
\Proof Again, Looijenga's proof of the corresponding result, 
\cite{looijenga} 8.20, applies practically verbatim. 
The proof is by induction on 
$\dim\ S$, but does not involve fibrations with fibre dimension
greater than $p$. There is a minor difference of notation: Looijenga
uses $\o_f$ in place of Greuel's $H '''$, to which our notation
$\chh '''$ refers implicitly.
\eop
\bc
\section{The Gauss-Manin connection}
\label{gm}
\ec
Let $\pi:D\to S$ be a good Stein representative of 
a deformation which frees the
$p$-dimensional almost free divisor $D_0,x_0$, and let ${\cal B}\subset S$ 
be the logarithmic discriminant.

The sheaf $\RR^p\pi_\ast\CC_D$ is a local system off $\cal B$, and so
$\RR^p\pi_\ast\CC_D\otimes\OO_S$ is naturally endowed with a flat connection.
Because of the natural isomorphism 
$$\RR^p\pi_\ast\CC_D\otimes\OO_S\simeq {\cal H}^p
(\pi_\ast\ch^p_{D/S})$$ off $\cal B$, ${\cal H}^p
(\pi_\ast\ch^p_{D/S})$ thus has a natural flat connection off $\cal B$, the 
topological connection. 
In \cite{looijenga}
the extension of the corresponding flat connection on 
the vanishing cohomology of the Milnor fibration $f:X\to S$ of an ICIS 
to the singular Gauss-Manin connection 
on $H^p(\O^\b_{X/S,(x_0,0)})$ is described in terms of Lie derivatives; 
the same construction serves in the case of an admissible deformation of an
almost free divisor. Following Looijenga's account, we now sketch this. 

First, if $s\notin {\cal B}$, then every vector field germ $\eta\in\t_{S,s}$
has a lift to a section $\xi$ of $\pi_\ast\Der(\log D)_{s}$; as in \cite
{looijenga}, this is a consequence of the fact that $\pi:D\to S$ is a Stein
map. For $[\o]\in {\cal H}^p(\pi_\ast\ch^\b_{D/S})_{s}$, set
$$\nabla_\eta([\o])=[L_\xi(\o)],$$
where $L_\xi$ is the Lie derivative.
Note that Cartan's formula $L_{\xi}(\o)=d\iota_\xi(\o)+\iota_\xi(d\o)$ 
holds for $\o\in \ch^p_D$, since both the exterior derivative
$d$ and contraction by the logarithmic vector field $\xi$ preserve the 
subcomplex
$h\O^\b(\log D)$ of $\O^\b_A$ (where $A$ is the ambient Milnor ball). It
is easy to check that $L_\xi$ is well defined on the relative
complex $\O^\b_{D/S}$ and
commutes with its exterior derivative.

If $\xi'$ is another lift of $\eta$ then $\xi-\xi'$ is vertical with respect
to $\pi$, so
since $d\o\in\pi_\ast(\pi^\ast\O^1_{S,s}\wedge\ch^p_D)_{s}$,
we have 
$$L_{\xi-\xi'}(\o)=\iota_{\xi-\xi'}(d\o)+d(\iota_{\xi-\xi'}(\o))$$
$$= d(\iota_{\xi-\xi'}(\o)).$$
Thus $\nabla:{\cal H}^p(\pi_\ast\ch^\b_{D/S})_{s}\to
{\cal H}^p(\pi_\ast\ch^\b_{D/S})_{s}\otimes_{\OO_S}\O^1_{S,s}$ 
is well defined.

Looijenga's proof that $\nabla$ so defined coincides with the topological
connection works verbatim: one needs only to check that the pairing
${\cal H}^p(\pi_\ast\ch^\b_{D/S,s})\ \times H_p(D_{s};\CC)\to \OO_{S,s}$
(given by integrating forms over the translation of cycles 
$Z\in D_s$ to nearby fibres by
means of the local trivialisation of the fibration $\pi:D\to S$ around
$s$) is non-degenerate in ${\cal H}^p(\pi_\ast\ch^\b_{D/S,s})$,
and the argument on pages 149-150 of \cite{looijenga} shows this,
without any modification.\\

\ni The connection $\nabla$ extends to a meromorphic connection
on ${\cal H}^p(\pi_\ast\ch^\b_{D/S,s})$ for points $s\in{\cal B}$, as
follows: let $g$ be a defining equation for the hypersurface ${\cal B}$, 
and let ${\cal L}$ denote the (coherent) 
$\OO_S$ module of $\pi$-liftable vector
fields.  If $\eta\in\t_{S,s}$ then for some positive power $m$, 
$g^m\eta$ is liftable, by virtue of the coherence of $\t_S/{\cal L}$. 
Let $\xi$ be a lift.
Moreover, by the coherence of $\chch'$ and ${\cal H}^p(\pi_\ast\ch^\b_{D/S})$,
and the fact that they
coincide off $\cal B$, some power of $g$ pushes $L_{\xi}(\o)$ into
$Z^p(\pi_\ast\ch^\b_{D/S,s})$. Thus we can define 
$$\nabla:{\cal H}^p(\pi_\ast\ch^\b_{D/S})
\to {\cal H}^p(\pi_\ast\ch^\b_{D/S})[g^{-1}]\otimes\O^1_S$$
by 
$$\nabla_\eta([\o])=g^{-m}[L_{\xi}(\o)]$$
where $\xi$ is a lift of $g^m\eta$ and $m$ is a sufficiently high power 
for the expression on the right
to make sense.

Finally, this connection extends by Leibniz's rule to a 
meromorphic connection on
${\cal H}^p(\pi_\ast\ch^\b_{D/S,s})$.

\begin{theorem} The Gauss-Manin connection just defined has a  regular
singularity along $\cal B$.
\end{theorem}
\Proof Once again, the proof for the classical case 
given in \cite{looijenga} applies verbatim.
\eop

\bc
\section{Calculating $\mu_E$}\label{calc}
\ec
\begin{lemma} If $D_0,x_0$ is a $p$-dimensional 
almost free divisor based on the free
divisor $E$, and $p <\min\{\mbox{wh}(E),\mbox{hn}(E)\}$, 
then $\mu_E(D_0,x_0) = \dim\ \ch^p_{D_0,x_0}/d\ch^{p-1}_{D_0,x_0}$.
\end{lemma}
\Proof Let $\pi:D\to S$ be a deformation which frees 
$D_0,x_0$ (such deformations exist, by the hypothesis that
$p < \mbox{hn}(E)$).
Then $\chh ' = \ch^p_{D/S, (x_0,0)}/d\ch^{p-1}_{D/S,(x_0,0)}$
is a free $\OO_{S,0}$-module of rank $\mu_E$, and 
$\chh '/m_{S,0}\chh '
=\ch^p_{D_0,x_0}/ d\ch^{p-1}_{D_0,x_0}$.
\eop

\begin{proposition}\label{earl} Suppose that $D_0$ is a $p$-dimensional
almost free divisor based
on the free divisor $E$, and that $\pi:D\to S$ is a deformation over
the 1-dimensional base $S$,
which frees 
$(D_0,x_0)$, with the property that $(D,(x_0,0))$ is also almost free. Then
provided $p+1 <\mbox{wh}(E)$,
$\mu_E(D_0)+\mu_E(D)=\dim\ \ch^{p+1}_{D/S}$
\end{proposition}
\Proof The proof is practically identical to Greuel's proof of 
the corresponding result (Lemma 5.3) of \cite{greuel}. 
Consider the map
$$\nabla'_{d/ds}:ds\wedge \check{H}'\to\check{H}'', \ \ 
[ds\wedge\o]\mapsto [d\o].$$
Both $ds\wedge \check{H}'$ and $\check{H}''$ are $\CC\{s\}$-modules
of rank $\mu_E(D_0)$, with $ds\wedge \check{H}' \subseteq \check{H}''$, 
and 
$$\nabla'_{d/ds}(g ds\wedge[\o])
=g.\nabla'_{d/ds}(ds\wedge[\o])+(dg/ds)\nabla'_{d/ds}(ds\wedge[\o])$$
for $g\in\CC^{s}$; it follows by the Malgrange index theorem (see
for example \cite{Malgrange},Theorem 2.3)
that 
$$\dim \ker \nabla'_{d/ds} - \dim\ \mbox{coker} \nabla'_{d/ds}
=\mbox{rank}\ ds\wedge \check{H}'-\dim\ (\check{H}''/ds\wedge \check{H}').$$
Now $\nabla'_{d/ds}$ is injective 
(by~\ref{exacta}) 
and so this equality amounts to
$$-\dim\ \frac{\ch^{p+1}_{D,x}}{d\ch^p_{D,x}}=\mu_E(D_0)+\dim\  
\ch^{p+1}_{D/S},$$
and thus to the required equality.
\eop

\begin{lemma}\label{dims} $\dim\ \ch^{p+1}_{D/S,x}=\dim T^{1,\log}_{D/S,x}$.
\end{lemma}
\Proof Let $\rho:{\cal D}\to S\x T$ 
be a free extension of $\pi:D\to S$. 
By the proof of~\ref{torsion}(4), 
$T\ch^{p+1}_{D,x}=\ker\ (\lambda:\ch^{p+1}_{D,x}\to
\ch^{p+d+1}_{{\cal D},x})/(t)\ch^{p+d+1}_{{\cal D},x}$ 
As $\lambda d\pi:\ch^{p}_{D/S,x}\to \ch^{p+d+1}
_{{\cal D},x}$ is injective, 
$T\ch^{p+1}_{D,x}\cap d\pi(\ch^{p}_{D/S,x})=0$, and hence
there is an exact sequence
$$0\to T\ch^{p+1}_{D,x}\to \ch^{p+1}_{D/S,x}\to \frac{dt\wedge\ch^{p+1}_{D,x}}
{dt\wedge d\pi\wedge\ch^p_{D/S}}\to 0.$$
From the exact sequence
$$0\to \frac{dt\wedge\ch^{p+1}_{D,x}}
{dt\wedge d\pi\wedge\ch^p_{D/S}}\to T^{1,\log}_{D/S}\to T^{1,\log}_{D}
\to 0$$ 
and the fact that $\dim\ T\ch^{p+1}_{D,x}=\dim T^{1,\log}_{D}$ 
(by \ref{torsion})
we conclude that
$\dim \ch^{p+1}_{D/S,x}=\dim T^{1,\log}_{D/S,x}$.
\eop
\begin{corollary}\label{count} With the hypotheses of \ref{earl},
 $\mu_E(D_0)+\mu_E(D)=\dim_{\CC}\ T^{1,\log}_{D/S}$.
\eop
\end{corollary}
Exactly as in \cite{greuel}, this furnishes us with a means of calculating
$\mu_{E}(D_0)$: we suppose that we are given a free 
deformation $\pi:D\to S$ of $D_0$ with the property that if 
$D_i=D\cap\{s_{i+1}=\cdots=s_d=0\}$ then each germ $(D_i,0)$ is almost free, 
and with $\dim D < \mbox{wh}(E)$. We
apply \ref{count} to the deformations $\pi_i:D_i\to S_i$, where $S_i$ is 
1-dimensional and $\pi_i(x,s)=s_i$. 
\begin{corollary}\label{duke} $$\mu_E(D_0,0)=
\sum_{i=1}^d(-1)^{i+1}\dim_{\CC}(\ch^{p+i}_{D_i/S_i,0})$$
$$=\sum_{i=1}^d(-1)^{i+1}\dim_{\CC}(T^{1,\log}_{D_i/S_i}).$$
\eop
\end{corollary}
This should be contrasted with the calculation of $\mu_E$ given in
\cite{dm}, which we now summarise. 
Let $\pi:(D, (x_0,0))\to (S,0)$ be a free
deformation of the almost free divisor $D_0,x_0
= i^{-1}(E)$, suppose that $\pi$ frees $D_0$, 
and suppose that we are given a defining
equation $h$ for $D$ in the smooth ambient space $A$, with the 
property that there exists a vector field $\chi\in\t_A$ such that $\chi.h=h$
(such an equation is called a "good  defining equation" in \cite{dm}). Let
$\Der(\log h)$ be the set of vector fields annihilating $h$ - i.e., tangent
to all of the level sets of $h$. Then we have
\begin{proposition}\label{prince} In these circumstances
$$\mu_E(D_0,x_0)=\dim_{\CC}\ \frac{\t(i)}{ti(\t_V)+i^\ast(\Der(\log h))}$$
$$=\dim_{\CC}\ \frac{\t(\pi)}{t\pi(\Der(\log h))+m_{S,0}\t(\pi)}.$$
\eop
\end{proposition}
In fact \cite{dm} proves only the first of these equalities. The 
second follows, by the argument of the proof of \ref{iso}.

Our \ref{count} and  \ref{duke} do not reduce in any obvious way to 
 \ref{prince}. For example, consider the case of a free and freeing
deformation $\pi:D\to S$ with $\dim S=1$. Then  \ref{count} gives
\begin{equation}\label{pip}
\mu_E(D_0,x_0)=\dim\ \frac{\t(\pi)}{t\pi(\Der(\log D))}
\end{equation}
while  \ref{prince} gives
\begin{equation}\label{pop}
\mu_E(D_0,x_0)=\dim\ \frac{\t(\pi)}{t\pi(\Der(\log h))+
m_{S,0}\t(\pi)}.
\end{equation}
Moreover,  \ref{count} does not require a good defining equation.

In fact the modules on the right hand side of (\ref{pip}) and (\ref{pop})
coincide when $D$ is weighted homogeneous and the variable $s$ has non-zero
weight $w_s$: if
$\chi_e=w_ss\p /\p s+\sum_i w_ix_i\p /\p x_i$ is the Euler vector field, then
$\Der(\log D)=\Der(\log h)\oplus \langle \chi_e\rangle$, and 
$t\pi(\Der(\log D))=t\pi(\Der(\log h))+m_{S,0}\t(\pi)$. 

Another application of \ref{count} is the following:

\begin{corollary} If $D_0$ is a $p$-dimensional almost free divisor
based on the free divisor $E$, if $p<\min\{\mbox{wh}(E),\mbox{hn}(E)\}$,
and if $\mu_E(D_0)=1$, then the ${\cal K}_E$-discriminant in a versal
deformation of $D_0$ is reduced.
\end{corollary}
\Proof Let $i_0$ induce $D_0$ from $E$.
By the main theorem of \cite{dm}, the ${\cal K}_{E,e}$-codimension of 
$i_0$ is no greater than $\mu_E(D_0)$; it is therefore equal to 1. Let
$i$ be a miniversal deformation of $i_0$ over the smooth 1-dimensional base
$S$, and let $D=i^{-1}(E)$.

The ${\cal K}_E$ discriminant consists simply of $\{0\}$, but with
analytic structure provided by the $0$-th Fitting ideal of $T^{1,\log}_{D/S}$
as $\OO_S$-module. By \ref{count}, $\mu_E(D_0)=\dim_{\CC}T^{1,\log}_{D/S}$
(since $i$ is algebraically transverse to $E$). Therefore $T^{1,\log}_{D/S}$
has length 1. It follows that its zero'th Fitting ideal over $\OO_S$ is
the maximal ideal $m_{S,0}$; hence the ${\cal K}_E$-discriminant is
reduced.\eop

We remark that the same argument shows that the vector field $s\partial/\partial s$ on $S$ is liftable to a vector field  $\chi\in\Der(\log D)$. See
\cite{damonlegacy} for some striking results on the question of the
reducedness and freeness of the ${\cal K}_E$-discriminant.

\bc
\section{Discriminants of Maps and Damon's theorem}\label{last}
\label{maps}
\ec
In this section we describe the relation between singularities and
deformations of map-germs $f:\CC^n,0\to\CC^p,0$ and 
the theory of almost free divisors. 

In \cite{damonKV}, Damon proves the following theorem:
\begin{theorem}\label{AKV} Suppose that the diagram 
$$\begin{array}{ccc}{\cal U},0&\stackrel{F}{\too}&{\cal V},0\\
j_0\uparrow\ \ &&i_0\uparrow\ \ \\
U,0&\stackrel{f}{\too}&V,0
\end{array}$$
is a fibre square (with $i_0$ transverse to $F$), and that $F$ is stable.
Let $D$ be the discriminant of $F$ (or image if $n<p$). Then 
$$T^1_f\simeq \frac{\t(i_0)}{ti_0(\t_{\CC^p})+i_0^\ast(\Der(\log D))}.$$
\end{theorem}
Here $T^1_f$ is the $\OO_{\CC^p,0}$-module
$$\frac{\t(f)}{tf(\t_{\CC^n})+\o f(\t_{\CC^p})},$$
also known as $N{\cal A}_ef$, and all the spaces are assumed non-singular.

Damon's proof in \cite{damonKV} used the additional hypothesis
that $f$ has finite
$\A_e$-codimension, and was rather involved. C.T.C.
Wall gave a simpler proof in \cite{walldim}; here we give a proof along
similar lines. Other proofs can be found in \cite{GdPW} and \cite{MM}.\\

\ni{\bf Proof of \ref{AKV}}\hspace{.3in} 
Let us first assume that $F$ is an unfolding of
$f$, with ${\cal U}=U\x S, {\cal V}=V\x S$, and that 
$i_0:V,0\to {\cal V},0$ is the
standard inclusion $v\mapsto (v,0)$.\\
{\bf Step 1}\ \ 
$$\frac{\t (F)}{tF(\t_{U\x S})}\ \simeq\ \frac{\t(F/\pi)}
{tF(\t_{U\x S/S})}.$$
{\bf Step 2} \  
By a standard argument (see e.g. \cite{looijenga} 6.14), $\Der(\log D)$
is the kernel of the morphism $\o F:\t_{V\x S}\to
F_\ast(\t(F)/tF(\t_{U\x S}))$. There is a commutative diagram
\renewcommand{\arraystretch}{1.5}
$$\begin{array}{ccccccccccc}
&&0&&0&&&&&&\\
&&\uparrow&&\uparrow&&&&&&\\
0&\to&t\pi(\Der(\log D))&\to&\t(\pi)&\to&\frac{\t(\pi)}{t\pi(\Der(\log D))}&
\to&0&&\\
&&t\pi\ \uparrow \ &&t\pi\ \uparrow &&&&&&\\
&&&&&&&&&&\\
0&\to&\Der(\log D)&\to&\t_{V\x S}&\to&F_\ast\left(\frac{\t(F)}
{tF(\t_{U\x S})}\right)&\to&0&&\\
&&\uparrow&&\uparrow&&\uparrow\simeq&&&&\\
0&\to&\Der(\log D)_\pi&\to&\t_{V\x S/S}&\to&F_\ast\left(\frac{\t_\pi(F)}
{tF(\t_{U\x S/S})}\right)&\to&T^1_{F/S}&\to&0\\
&&\uparrow&&\uparrow&&&&&&\\
&&0&&0&&&&&&
\end{array}$$ 
\renewcommand{\arraystretch}{1}
Diagram chasing gives an isomorphism
$$\frac{\t(\pi)}{t\pi(\Der(\log D))}\simeq T^1_{F/S};$$
so now tensoring with $\OO_{V,0}$ (over $\OO_{V\x S,0}$) produces
an isomorphism
$$\frac{\t(\pi)}{t\pi(\Der(\log D))+m_{S}\t(\pi)}=
\frac{\t(\pi)}{t\pi(\Der(\log D))}\otimes\OO_{V,0}\simeq T^1_f.$$

But
$$\frac{\t(\pi)}{\pi^\ast m_{S,0}\t(\pi)}$$
can be identified with
$$\frac{\t(i_0)}{ti_0(\t_V)}$$
and using this identification we obtain
$$\frac{\t(i_0)}{ti_0(\t_V)+i_0^\ast\Der(\log D)}\simeq\frac{\t(\pi)}
{t\pi(\Der(\log D))+m_{S,0}\t(\pi)}\simeq T^1_f.$$

This completes the proof in the special case where $F$ is a (parametrised)
unfolding of $f$. The general case follows by showing first 
that up to isomorphism,
the module $$\frac{\t(i)}{ti(\t_V)+i_0^\ast\Der(\log D)}$$ is independent
of the choice of stable map $F$ from which $f$ is pulled back by $i$. For
stable maps are classified by their local algebras (\cite{mather4}),
and the local algebra of $F$ is isomorphic to that of $f$;  hence, if 
$F_1:{\cal U}_1,0\to {\cal V}_1,,0$ and 
and $F_2:{\cal U}_2,0\to {\cal V}_2,0$
are  stable maps from which $f$ may be induced by transverse 
pull-back, with $\dim\ {\cal U}_1\geq \dim\ {\cal U}_2$, 
then $F_1$ is equivalent to a trivial
unfolding of $F_2$ on $\dim\ {\cal U}_1-\dim\ {\cal U}_2$ parameters.
As a consequence of this, we may assume that $F$ is a stable unfolding
of $f$, as in the special case.\eop

Note that the isomorphism $T^1_{F/S}\simeq \t(\pi)/t\pi(\Der(\log D))$
identifies the support of $T^1_{F/S}$ 
as the logarithmic critical space of $\pi:D\to S$.
Note also that there is now a strong formal analogy between $T^1_f$ and
$T^1_{X_0}$, where $X_0=f^{-1}(0)$:
$$T^1f=\frac{\t(\pi)}{t\pi(\Der(\log D))+ m_{\CC^d,0}\t(\pi)}$$
$$T^1X_0=\frac{\t(f)}{tf(\t_{\CC^n,0})+m_{\CC^p,0}\t(f)}.$$

Author's address:\\
Mathematics Institute, University of Warwick, Coventry CV4 7AL, United Kingdom\\
e-mail: mond@maths.warwick.ac.uk
\end{document}